\documentclass[preprint,12pt]{elsarticleplain}

\usepackage{caption}
\usepackage{subfigure}
\usepackage{amssymb}
\usepackage{amsthm}
\usepackage{amsmath}
\usepackage{mathtools}
\usepackage{algorithm}
\usepackage[noend]{algpseudocode}
\usepackage[hidelinks]{hyperref} 
\usepackage{diagbox}
\usepackage{multirow}

\newtheorem{theorem}{Theorem}[section]
\newtheorem{lemma}{Lemma}[section]
\newtheorem{proposition}{Proposition}[section]
\newtheorem{remark}{Remark}[section]
\newtheorem{definition}{Definition}[section]
\newtheorem{assumption}{Assumption}[section]
\newtheorem{example}{Example}[section]

\numberwithin{equation}{section}

\def\sin{\mathrm{sin}}
\def\cos{\mathrm{cos}}
\def\e{\mathrm{e}}

\usepackage{xcolor}
\newcommand{\RomanNumeralCaps}[1]{\MakeUppercase{\romannumeral #1}}

\DeclarePairedDelimiter\ceil{\lceil}{\rceil}
\DeclarePairedDelimiter\floor{\lfloor}{\rfloor}
\DeclareMathOperator*{\argmin}{arg\,min}

\makeatletter
\def\@author#1{\g@addto@macro\elsauthors{\normalsize%
    \def\baselinestretch{1}%
    \upshape\authorsep#1\unskip\textsuperscript{%
      \ifx\@fnmark\@empty\else\unskip\sep\@fnmark\let\sep=,\fi
      \ifx\@corref\@empty\else\unskip\sep\@corref\let\sep=,\fi
      }%
    \def\authorsep{\unskip,\space}%
    \global\let\@fnmark\@empty
    \global\let\@corref\@empty  
    \global\let\sep\@empty}%
    \@eadauthor={#1}
}
\makeatother
\begin{document}

\begin{frontmatter}



\title{A fast neural hybrid Newton solver adapted to implicit methods for nonlinear dynamics}


\author[inst1,inst3]{Tianyu Jin\corref{cor1}}
\ead{tjinac@connect.ust.hk}
\cortext[cor1]{Corresponding author}
\author[inst2]{Georg Maierhofer}
\author[inst3]{Katharina Schratz}
\author[inst1,inst4]{Yang Xiang}

\affiliation[inst1]{organization={Department of Mathematics},
            addressline={The Hong Kong University of Science and Technology}, 
            state={Clear Water Bay},
            country={Hong Kong Special Administrative Region of China}}

\affiliation[inst2]{organization={Mathematical Institute},
            addressline={University of Oxford},
            country={United Kingdom}}

\affiliation[inst3]{organization={Laboratoire Jacques-Louis Lions (UMR 7598)},
            addressline={Sorbonne Université},
            country={France}}

\affiliation[inst4]{organization={Algorithms of Machine Learning and Autonomous Driving Research Lab},
            addressline={HKUST Shenzhen-Hong Kong Collaborative Innovation Research Institute}, 
            city={Futian},
            state={Shenzhen},
            country={China}}

\begin{abstract}
The use of implicit time-stepping schemes for the numerical approximation of solutions to stiff nonlinear time-evolution equations brings well-known advantages including, typically, better stability behaviour and corresponding support of larger time steps, and better structure preservation properties. However, this comes at the price of having to solve a nonlinear equation at every time step of the numerical scheme. In this work, we propose a novel deep learning based hybrid Newton’s method to accelerate this solution of the nonlinear time step system for stiff time-evolution nonlinear equations. We propose a targeted learning strategy which facilitates robust unsupervised learning in an offline phase and provides a highly efficient initialisation for the Newton iteration leading to consistent acceleration of Newton's method. A quantifiable rate of improvement in Newton's method achieved by improved initialisation is provided and we analyse the upper bound of the generalisation error of our unsupervised learning strategy. These theoretical results are supported by extensive numerical results, demonstrating the efficiency of our proposed neural hybrid solver both in one- and two-dimensional cases. 
\end{abstract}




\end{frontmatter}


\section{Introduction}\label{sec:intro}
We consider the numerical approximation of solutions to stiff nonlinear time-evolution partial differential equations (PDEs) with implicit numerical methods. In particular, the type of equation we consider in this manuscript is of the following general form with appropriate boundary conditions (e.g. periodic, homogeneous Neumann or Dirichlet)
\begin{equation}\label{eq...evol}
\left\{\begin{array}{l}
\partial_tu=Au+f(u), \textbf{x}\in\Omega, t\in[0,T],\\
u(\textbf{x},0) = u_0(\textbf{x}), \textbf{x}\in\overline{\Omega},\\
\end{array}\right.
\end{equation}
where $\Omega\subset\mathbb{R}^d(d = 1,2,3)$ is an open and bounded domain, $u:\overline{\Omega}\times[0,T]\rightarrow\mathbb{R}$ is the unknown function, $A$ is a stiff spatial differential operator which causes the stiffness and $f:\mathbb{R}\rightarrow\mathbb{R}$ is a nonlinear function. For our purposes `stiff' in this context will refer to high order spatial derivative operators that inherently require implicit time-steppers to overcome stability issues \cite{hairer1993solving}.

While the stiff nature of \eqref{eq...evol} makes this problem generally challenging, over the last few decades a wide range of successful numerical methods for the approximate solution of \eqref{eq...evol} have been developed. Explicit methods, while often simpler and more computationally efficient, are typically conditionally stable. For example, a necessary condition, now known as the CFL condition \cite{CFL}, needs to be satisfied for the convergence of the explicit finite difference schemes. On the other hand, implicit methods, and in particular fully implicit methods, are generally unconditionally stable. This implies that they are usually more stable and have better structure-preserving property with wider choice of time step or space increment \cite{hairer1993solving}. Nevertheless, some iterative methods are required in solving fully implicit methods at each time step, which makes them more computationally expensive than explicit methods.

A multitude of iterative methods exist for the resolution of systems of nonlinear equations, each with their unique characteristics and applications. Fixed-point iteration, for instance, is the simplest one while the strong convergence requirement restricts its usage. Picard iteration is another commonly used approach by successively approximating the solution as an integral involving the unknown function. Among all the iterative methods, Newton iteration is a powerful technique with robustness and quadratic convergence rate. However, the evaluation of the Jacobian matrix is indispensable within each iteration of Newton's method, thereby causes substantial computational cost. Furthermore, it is noteworthy that Newton's method exhibits high sensitivity to the initial guess (see Section \ref{sec:estimates_initialisation_in_newtons_method} for detailed explanation), as an initial guess that is far from the true solution may result in failure to converge. Hence, there arise two key challenges for improving Newton's method: (1) mitigating the computational cost, and (2) devising strategies to obtain a good initial guess. Several conventional approaches have been developed to address the first challenge. Quasi-Newton methods \cite{quasinewton,quasireview} are a class of methods based on establishing an approximation to the Jacobian matrix using successive iterates and function values. However, compared to Newton’s method, quasi-Newton methods generally require a larger iteration count. Thus, the cost of these extra iterations is usually higher than the savings in approximating Jacobian per iteration. Another comparable class of methods are so-called inexact Newton methods, originally proposed in \cite{inexact}, which do not require the exact computation of each Newton iterate and instead allow for a small error at each iteration (for example caused by iterative solvers for the corresponding linear Newton-system) by introducing a relative residual into the Newton equation. A notable limitation of inexact Newton methods lies in the constant trade-off between the accuracy needed to maintain Newton iteration's convergence rate and the computational cost per iteration. This implies that although a large tolerance of error due to residuals can decrease the computational cost per iteration, it may lead to an increase of total iteration count because of the lower convergence rate. Later on, other methods are also proposed based on these two classes of methods. One example is the Newton-Krylov-Schwarz solver \cite{LIU2024112548} for the coupled Allen--Cahn/Cahn--Hilliard system, which combines the inexact Newton iteration with an additive Schwarz preconditioner to enhance computational efficiency. However, choosing proper boundary conditions for the subdomain problems has a huge impact on the convergence of the Schwarz preconditioner. Hence, this method is problem-specific and has some difficulty in dealing with other problems.

In the last decade, deep learning has gained considerable success across various areas of science and engineering, especially in solving PDEs. One notable advancement in this area is operator learning \cite{chenchen,Lulu,graphkernelnet,li2021fourier}, i.e., using the neural network to approximate the mappings between two infinite-dimensional function spaces. Although operator learning exhibits better generalisation property compared to directly learning the unknown function, it still suffers from the problem of non-guaranteed convergence for approximated solutions, unlike conventional numerical methods. As the time step and the spatial increment tend to $0$, the error of the neural network approximated solution does not converge to $0$. As a result, it would be more advantageous to combine deep learning methods with classical numerical methods in order to maximize the benefits of both. The use of deep learning techniques for improving iterative PDE solvers has been documented in several previous works, see \cite{Solver-in-the-Loop, HINTS}. \cite{qianxiao1, qianxiao2} also consider utilizing meta-learning approaches to speed up iterative algorithms, such as the Jacobi method and multigrid method. When it comes to Newton's method, \cite{IntDeep} introduced a two-phase algorithm for solving elliptic equations. In phase 1, a deep learning solver is adopted to solve the equations roughly. While in phase 2, the network prediction is converted to a finite element ansatz, which then serves as an initial guess for Newton's method to solve the weak form of nonlinear equations.
 \cite{doi:10.1137/22M1507942} proposed a nonlinearly preconditioned inexact Newton method with learning capability where a PCA-based unsupervised learning strategy is employed to learn a decomposition in constructing preconditioner. It would also be possible to address the second challenge of Newton's method by using neural operators due to their acceptable accuracy regardless of the time step. In \cite{co2,aghili2024accelerating}, Fourier Neural Operator (FNO) is trained as a time stepper to give a better initial guess in Newton's method for solving different nonlinear systems of equations. However, training FNO requires a large amount of labeled data, which is often expensive to acquire. Furthermore, inference over FNO typically demands more CPU time compared to other neural network structures due to the inclusion of the Fast Fourier Transform (FFT) operator. Therefore, the challenge of efficiently training a neural network to furnish a better initial guess in Newton iteration when dealing with time-evolution equations emerges as a significant problem.

In this work, we propose a general neural hybrid solver for a class of stiff nonlinear time-evolution equations. Unlike the supervised learning strategy used in \cite{co2,aghili2024accelerating}, an implicit-scheme-informed unsupervised learning method is constructed to train our neural time stepper with extremely simple and light structure compared to most other neural PDE solvers. Then the output of the neural network is used as the initial guess in Newton's method, which can decrease the total iteration count and accelerate the whole algorithm in solving \eqref{eq...evol} with relatively large time step. Notice that once the neural time stepper is trained, it can consistently provide a good initial guess in multiple time steps until equilibrium without being trained again. We evaluate the proposed method on the Allen--Cahn equation both in $1$- and $2$-dimensional cases. In such problems, fully implicit methods usually have better structure-preserving properties than explicit methods but are generally more computationally expensive. Examples of several state-of-the-art implicit methods for this equation which are energy-stable include, \cite{tierra_review_2015,timestepadapative,gomez2011,hou2023energy}. Based on the implicit midpoint method, \cite{AC_im_midpoint} proposed an implicit step-truncation midpoint method and experimental results indicate its efficiency in solving Allen-Cahn equation. Numerical results demonstrate that the proposed neural hybrid solver outperforms the classical Newton solvers in terms of computational efficiency. While the hybrid approach taken in this work guarantees that the final output of Newton iterations (provided it converged) inherits the rigorously guaranteed properties of the original implicit scheme, we provide further theoretical bounds on the behaviour of the hybrid Newton method. 
In particular, we provide quantitative estimates on the effect of an improved initialisation on the iteration count in Newton's method, results which are found to closely match the behaviour in practice.
Moreover, we estimate the generalisation error of our proposed unsupervised learning strategy which can be bounded by the training error and the covering error of the training dataset on the underlying function space. This estimate is comparable to the approach taken recently in \cite{mishra22} but differs since in our case we quantify the generalisation error of an operator while \cite{mishra22} considered essentially a function approximator. 

The rest of the paper is structured as follows: in Section \ref{sec:hybrid_methods}, we introduce the proposed neural hybrid solver in detail. Section \ref{sec:theoretical} presents the asymptotic iteration count estimation of Newton’s method with exponential decaying initialisation and a theoretical analysis of the generalisation error of the employed unsupervised learning strategy under appropriate assumptions. The numerical implements on solving $1$- and $2$-dimensional Allen--Cahn equations are shown in Section \ref{sec:numerical_examples}. 

\section{Neural hybrid Newton solver}\label{sec:hybrid_methods}
In this section, we introduce the proposed neural hybrid method to solve the implicit time-stepper. In order to solve \eqref{eq...evol}, we first divide the time interval $[0,T]$ into sub-intervals $0=t_0<t_1<\ldots<t_{N_T}=T$ such that $t_n-t_{n-1}=\tau$ for all $n=1,2,\dots,N_{T}$. Suppose $u^n\in \mathbb{R}^{N^d}$, $u^{n+1}\in \mathbb{R}^{N^d}$ are the approximate solution values at $t_n$, $t_{n+1}$ under the $N^d$ uniform spatial discretisation. Then we denote the implicit method concerned solving \eqref{eq...evol} to be
\begin{equation}\label{scheme...im}
u^{n+1} = \Psi_{\tau}(u^n,u^{n+1})
\end{equation}
and we further assume the implicity is involved in the nonlinear term in the operator $\Psi_{\tau}$ (i.e. fully implicit scheme). Therefore, solving the system of nonlinear equations $u=\Psi_\tau(u^n,u)$ is required at each time step. A natural approach is to use Newton iteration. We start our description by recalling the precise definition of Newton iteration.
\begin{remark}\label{rmk...Psi}
    To denote a slight abuse of notation, we use $\Psi_\tau$ to refer both to its discretised form and to the operator applied to variables considered as continuous functions. This dual usage should be clear from the context in which $\Psi_\tau$ is applied. The similar dual usage is applied to $u^n$. As will become apparent throughout the following sections, our theoretical analysis is based on functions which are then discretised in the numerical implementation of our hybrid method. Whenever appropriate the reader may wish to interpret this dual notation through the lens of interpolation versus discretisation - in our numerical examples with homogeneous Neumann boundary conditions we use a central finite difference scheme in space and interpolation over the Lagrange finite element space of degree two, while for periodic boundary conditions one may think of discretisation and interpolation based on discrete Fourier transforms and trigonometric polynomials.
\end{remark}

\subsection{Newton iteration}\label{sec:def_newton}
We define a new operator $G(\textbf{y}):=\textbf{y}-\Psi_{\tau}(u^n,\textbf{y})$ for all $\textbf{y}\in \mathbb{R}^{N^d}$ mapping from $\mathbb{R}^{N^d}$ to $\mathbb{R}^{N^d}$. Then solving the semi-discretised version of \eqref{scheme...im} is equivalent to solving 
$$
G(\textbf{y})=\textbf{0}, \quad\textbf{y}\in \mathbb{R}^{N^d}.
$$
The Newton iteration of solving systems of nonlinear equations $G(\textbf{y})=\textbf{0}$ is defined as follows.
\begin{definition}\label{def...newton}
Denote $D_\textbf{y}G(\textbf{y})$ the Jacobian matrix of $G$ and we further assume $D_\textbf{y}G(\textbf{y}^m)$ exists and is non-singular for each $m=0,1,2,\dots$. Then the recursion defined by
$$
\textbf{y}^{m+1}=\textbf{y}^m - D_\textbf{y}G(\textbf{y}^m)^{-1}G(\textbf{y}^m), \quad m = 0,1,2,\dots
$$
where $\textbf{y}^0\in \mathbb{R}^{N^d}$, is called Newton's method (or Newton iteration) for the system of equations $G(\textbf{y})=\textbf{0}$.
\end{definition}
As will become a bit more apparent in Section~\ref{sec:estimates_initialisation_in_newtons_method} the choice of the initial guess $\textbf{y}^0$ is crucial for the success and convergence speed of Newton's method. A direct way is to use $u^n$ as the initial guess. However, if $\tau$ is large in \eqref{scheme...im}, the direct initial guess $u^n$ may be far from the final solution $u^{n+1}$ and thus cause a large iteration count or, even worse, prevent Newton's method from converging. In the following, we establish a novel neural time stepper approximating the operator mapping from $u^n$ to $u^{n+1}$ to obtain a better initial guess in Newton's method.

\subsection{Neural network architecture}\label{sec:nn_architect}
In this subsection, we describe the architecture of the neural network in our method. The network we used is a simplified version of the commonly used convolutional neural network (CNN). Generally speaking, the whole neural network consists of $L$ convolutional layers, with a convolution operator and a nonlinear activation function inside each of them. Also, the number of channels $N_c$ of each convolution operator is increasing within the first $\floor*{\frac{L}{2}}$ layers and decreasing within the last $\floor*{\frac{L}{2}}$ layers, which gives us a structure that is wide in the middle and narrow at the ends. 
\begin{remark}
In the most commonly used CNNs, convolutional layers are usually followed by certain downsampling layers, such as Max Pooling and Average Pooling, in order to reduce the spatial dimensions of feature maps while retaining essential information. Since such downsampling layers are not included in our neural network structure, we call it a simplified CNN.
\end{remark}

Before giving the rigorous expression of the network architecture, we first review some basic operations in neural networks. We denote by $X\in \mathbb{R}^{m\times n}$ the input feature of the convolution and by $W\in \mathbb{R}^{k\times l}$ the weight of the convolution with $k\leq m$ and $l\leq n$. We use the notation $W* X$ to represent the correlation of $W$ and $X$, which reads as
$$
(W* X)_{i,j} = \sum_{r=1}^k\sum_{s=1}^l W_{r,s}X_{i+r-1,j+s-1}, \quad 1\leq i\leq m-k+1, 1\leq j\leq n-l+1.
$$

Then the notation $pad(X,p)\in \mathbb{R}^{(m+2p)\times(n+2p)}$ is used to represent the padding operator with padding size $p$. In our settings, we mainly consider homogeneous Neumann boundary condition. According to the conclusion in \cite{BoundaryConditionsCNN}, we choose the \textit{reflect} padding mode in the convolution operator which has the following form.
$$
(pad(X,p))_{i,j}=\left\{\begin{array}{l}
X_{i,j}\hspace{5em} p+1\leq i\leq p+m, p+1\leq j\leq p+n\\
X_{2+p-i,2+p-j}\hspace{1em} 1\leq i,j\leq p, m+p+1\leq i\leq m+2p,\\
                             \hspace{7em}m+p+1\leq j\leq n+2p\\
\end{array}\right.
$$

\begin{figure}[htbp]
	\centering
	\includegraphics[width=0.9\textwidth]{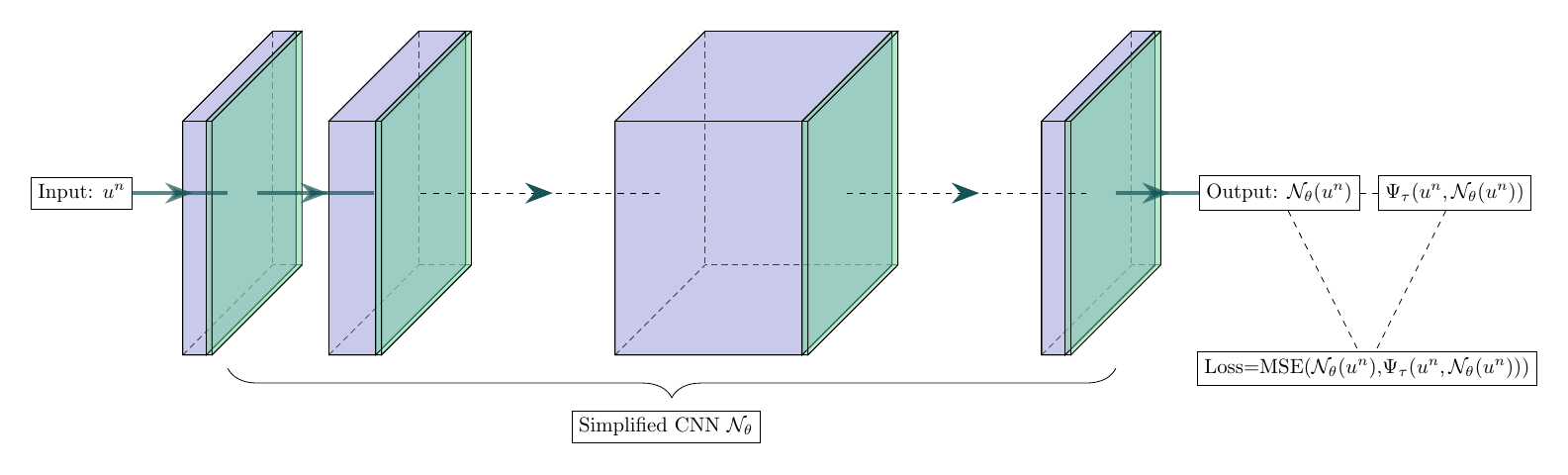}
	\caption{The architecture of the proposed implicit-scheme informed neural network. The purple blocks represent the convolution operator and the green blocks represent the activation function.}
	\label{fig:NN_architecture}
\end{figure}

Finally, let $\sigma$ denote the nonlinear activation function and $X^{(q)}$ denote the output of the $q$-th convolutional layer where $X^{(q)}_z\in \mathbb{R}^{m\times n}$ denote the $z$-th channel of it. Suppose $X^{(q-1)}$ has $N_{c1}$ channels and $X^{(q)}$ has $N_{c2}$ channels. Note that the number of channels should be $1$, if $q=1\text{ or }L$. Then the structure of each convolutional layer can be represented as follows.
\begin{align*}
    X^{(q)}_z = \sigma(\sum_{j = 1}^{N_{c1}}W_{j,z}^{(q)}* pad(X^{(q-1)}_j,p)\oplus b_z^{(q)}),
\end{align*}
with the trainable weights $W_{j,z}^{(q)}\in\mathbb{R}^{k\times l}$, biases $b_z^{(q)}\in\mathbb{R}$ and $\oplus$ denotes component-wise addition. The whole architecture with $L$ convolutional layers is illustrated in Figure \ref{fig:NN_architecture}.

An appropriate choice of neural network architecture is crucial and should be tailored to the specific task at hand. We hereby present an in-depth discussion on the benefits of employing such a simplified convolutional neural network structure for our task. As mentioned in the survey \cite{CNNreview}, a distinguishing feature of the convolution operation is the weight-sharing technique, whereby multiple sets of features within a data are extracted by sliding a kernel with the same set of weights. Therefore, compared to the fully-connected neural networks, the convolutional neural network possess two key advantages. Firstly, the reduced number of parameters in the network facilitates a faster and more efficient forward pass, as well as easier training of the network. Secondly, the network is able to detect patterns regardless of the location within the input. Moreover, the small spatial dimensions of the convolutional kernels enable the network to focus on local regions of the input and capture the correlations between neighbouring elements. This locality property has an inherent connection to the finite difference numerical schemes and is very well-suited to represent local differential operators present in the equations. We found in initial tests of different network architectures that our choice leads to a favourable balance between fast inference and accuracy for the problem at hand. (see Section \ref{sec:impact} for details.)

\subsection{Implicit-scheme-informed learning: Loss function and optimisation}\label{sec:scheme_informed_learning}
In this subsection, we introduce the training strategy of the aforementioned simplified CNN denoted by $\mathcal{N}_{\theta}$. Our target is to find an initial guess as close to $u^{n+1}$ as possible. Thus, we wish the network $\mathcal{N}_{\theta}$ can approximate one step of the implicit scheme \eqref{scheme...im}, i.e. approximate the operator $\Phi_\tau:u^n\mapsto u^{n+1}=\Phi_\tau(u^n)$ for all $n=0,1,2,...$. This is achieved by an implicit-scheme-informed unsupervised learning strategy as shown in Figure \ref{fig:NN_architecture}. An abstract paraphrasing of this strategy runs as follows. We define the residual:
\begin{equation}\label{eq...res}
    \mathcal{R}(u;u^n):= \Psi_{\tau}(u^n,u)-u
\end{equation}
Note that for solution $u^{n+1}$ of \eqref{scheme...im}, $\mathcal{R}(u^{n+1};u^n)=\Psi_{\tau}(u^n,u^{n+1})-u^{n+1}\equiv0$. Assume the input of the network is $u^n$ and $u^n, u^{n+1}\in X$ where $X$ is a function space. The training strategy is to minimize the residual \eqref{eq...res} over admissible set of trainable parameters $\theta\in\Theta$, i.e.
\begin{equation}\label{eq...train}
    \text{Find } \theta^*\in\Theta: \quad\theta^*=\argmin_{\theta\in\Theta}\|\mathcal{R}(\mathcal{N}_{\theta}(u^n);u^n)\|_X
\end{equation}
Since it is impossible to compute the norm $\|\cdot\|_X$ exactly, we need to approximate it by discretisation. For notational simplicity, we only write the step $u_0 \in X\rightarrow u^1\in X$ in the loss function without loss of generality. Given a set of initial data $\{u_0^{(i)}\}_{i=1}^{N_{data}}$ as the "quadrature points" in function space $X$, we determine the parameters $\theta\in\Theta$ via minimising the following loss function as the discrete version of \eqref{eq...train}:
\begin{equation}\label{eq...loss}
    \theta^* = \argmin_{\theta\in\Theta}\mathcal{L}(\theta):=\argmin_{\theta\in\Theta}\frac{1}{N_{data}}\sum_{i = 1}^{N_{data}}\| \Psi_{\tau}(u_0^{(i)},\mathcal{N}_{\theta}(u_0^{(i)}))-\mathcal{N}_{\theta}(u_0^{(i)})\|^2.
\end{equation}

\begin{example}
    If we work with periodic boundary conditions on the domain $\Omega=[0,1]^d$ in \eqref{eq...evol}, a sensible choice of $(X,\|\,\cdot\,\|)$ would be periodic Sobolev spaces with the natural discretisation of the norm arising from a pseudospectral discretisation in space.
\end{example}

In our experiments, the open-source machine learning library PyTorch \cite{pytorch} is used to construct and train the neural network. Following common practice we solve the minimisation problem \eqref{eq...loss} approximately by stochastic gradient descent, specifically using ADAM \cite{adam}. 

After the neural network is trained, we use the model's prediction, $\mathcal{N}_{\theta}(u^n)$, as the initial guess in Newton's method. Combined with the recursion of Newton iteration in Definition \ref{def...newton}, we summarise our neural hybrid solver in Algorithm \ref{alg...hybrid}.

\begin{algorithm}
\caption{Neural hybrid solver for $G(\textbf{y}) =0$ (one time step $u^n\rightarrow u^{n+1}$)}\label{alg...hybrid}
\begin{algorithmic}
\State $\textbf{y}=\mathcal{N}_\theta(u^n)$ \Comment{$\mathcal{N}_\theta$ provides the initial guess}
\State \textbf{User-specified input:} tolerance $\epsilon_{tol}$, and the maximum number of iterations allowed $N_{maxiter}$.
\State $l_{update}=2\epsilon_{tol}, k=0$
\While{$l_{update}\geq\epsilon_{tol}$ and $k<N_{maxiter}$}
    \State calculate $D_\textbf{y}G(\textbf{y})$, $\textbf{b}=G(\textbf{y})$  
    \State solve $\Delta \textbf{y}$ from $D_\textbf{y}G(\textbf{y})\Delta\textbf{y} = \textbf{b}$
    \State $\textbf{y} = \textbf{y}+ \Delta \textbf{y}$
    \State $l_{update} = \|\Delta\textbf{y}\|_{L^2}$
    \State $k=k+1$
\EndWhile
\State $u^{n+1} = \textbf{y}$
\end{algorithmic}
\end{algorithm}

\begin{remark}\label{rmk...alg}
    In Algorithm \ref{alg...hybrid}, the initial guess can be replaced by other choices, such as $u^n$ or approximated solution from other numerical methods, to obtain other types of Newton solvers. According to the experimental results presented in Section $4$, the neural hybrid solver achieves the best performance among all Newton solvers tested in terms of the trade-off between the computational cost of calculating initial guess and the reduction in iteration count and, hence, reduction in overall cost of the time-stepper. The GMRES algorithm can be applied to solve $D_\textbf{y}G(\textbf{y}_{old})\Delta\textbf{y} = \textbf{b}$ after wrapping the matrix-vector multiplication $D_\textbf{y}G(\textbf{y}_{old})\textbf{z}, \forall \textbf{z}\in\mathbb{R}^{N^d}$ into a \verb|LinearOperator|, the common interface for performing matrix-vector products in the open-source Python library SciPy.
\end{remark}

\section{Theoretical results}\label{sec:theoretical}
In this section we describe two central theoretical results underpinning our approach for the acceleration of Newton iterations. The first of those results describes specifically the speed-up of the Newton iterations that can be achieved with an improved initialisation. The second part provides bounds on the generalisation error of the neural network in different settings, thus ensuring that the accuracy we measure during training extends to relevant unseen data (solution values). Finally, we show the proofs of the main theorems given in the second part in the last subsection.
\subsection{Theoretical iteration count reduction in Newton's method through initialisation}\label{sec:estimates_initialisation_in_newtons_method}
In this part, we provide theoretical results on how much the iteration count can be reduced in Newton's method (defined in Definition \ref{def...newton}) with better initialisation. To begin with, we state some notations and assumptions that will be used later in this subsection. We recall that the system of equations we aim to solve is $G(\textbf{y})=\textbf{0}, \textbf{y}\in\mathbb{R}^{N}$ (the general case is 
 $\textbf{y}\in\mathbb{R}^{N^d}$, but here for notational simplicity we consider $d=1$ and note that of course the results extend verbatim to the case of $d>1$).
 Our assumptions on $G$ are standard in this context (cf. \cite{SuliMayers2003}) and are summarised in Assumption~\ref{assump...G}.
 \begin{assumption}\label{assump...G}
 The system of equations $G(\textbf{y})=\textbf{0}, \textbf{y}\in\mathbb{R}^{N}$ and our sequence of Newton iterates $\textbf{y}^{(k)}, k=0,\dots$ satisfies the following three assumptions:
 
\noindent
(i) The system of equations has a unique solution denoted by $\boldsymbol\xi$, i.e. $G(\boldsymbol\xi)=\textbf{0}$. Furthermore, there is an open neighbourhood of $\boldsymbol\xi$ such that in $N_{neighbour}(\boldsymbol\xi)$, $G$ is defined and continuous, all the second-order partial derivatives of $G$ are defined and continuous, and $D_{\textbf{y}}G(\boldsymbol\xi)$ is nonsingular. 

 \noindent
 (ii) The initialisation of our Newton iterates is in this neighbourhood, $\mathbf{y}^{(0)}\in N_{neighbour}(\boldsymbol\xi)$.

 \noindent
(iii)The Jacobian matrix $D_{\textbf{y}}G$ of $G$ exists and is nonsingular for each element of the sequence $\textbf{y}^{(k)}, k=0,\dots$.
 \end{assumption}

In the subsequent theorem, we recall the well-known result that Newton's method (as given in Definition \ref{def...newton}) achieves quadratic convergence rate subject to the aforementioned assumptions.
\begin{theorem}[{Theorem $4.4$ in \cite{SuliMayers2003}}]\label{thm...quadconv}
Suppose the system of equations $G(\textbf{y})=\textbf{0}, \textbf{y}\in\mathbb{R}^{N}$ and $N_{neighbour}(\boldsymbol\xi)$ satisfy the Assumption \ref{assump...G}. Denote the error of the $k$-th iteration $\varepsilon_k := \left\|\boldsymbol{y}^{(k)}-\boldsymbol{\xi}\right\|_{\infty}$. Then, 
\begin{align*}
        \varepsilon_{k+1} \leq C \varepsilon_k^2,
\end{align*}
where 
\begin{align}\label{eq:C}
        C=\frac{1}{2} N^2 A_G B_G,
\end{align}
$$
A_G=\max _{1 \leq i, j, l \leq N} \max _{\boldsymbol{y} \in \bar{B}_{\epsilon}(\boldsymbol{\xi})}\left|\frac{\partial^2 G_i}{\partial y_j \partial y_l}(\boldsymbol{y})\right|,\quad
 B_G=\max _{\boldsymbol{y} \in \bar{B}_{\epsilon}(\boldsymbol{\xi})}\left\|\left[D_{\boldsymbol{y}} G(\boldsymbol{y})\right]^{-1}\right\|_{\infty},
$$
for some $\epsilon>0$ such that $\bar{B}_{\epsilon}(\boldsymbol{\xi})\subset N_{neighbour}(\boldsymbol\xi)$.
\end{theorem}

The proof of Theorem \ref{thm...quadconv} can be found in \cite{SuliMayers2003}. Using this convergence rate, we can now investigate the effect of improving our initial guess  $\varepsilon_{0}$ by a factor of $2^{-n},n\in\mathbb{N}$ on the iteration count upper bound for convergence in Newton's method. 
\begin{theorem}[Theoretical upper bound of iteration count]\label{thm...asym_iteration_count}
    Suppose we are able to improve the initial estimate $\varepsilon_0$ by a factor of $2^{-n}$ where $n\in\mathbb{N}$. Then an upper bound of the iteration count, $M_{iter}=M_{iter}(n)$, required for Newton's method to achieve a fixed tolerance $\varepsilon_{tol}$, is given as follows.
    \begin{align}
        M_{iter}=\max\{0, \ceil*{C_1-\log_2|C_2-n|}\}.
    \end{align}
    The constants $C_1$ and $C_2$ depend on the properties of the first- and second-order partial derivatives of $G$ (i.e. the norm of inverse Jacobian and the norm of Hessian as outlined above in Theorem \ref{thm...quadconv}). Furthermore, $C_1$ depends on the desired tolerance and $C_2$ depends on the initial guess, but neither is affected by $n$. Note that we need $n>C_2$ for the convergence of the Newton iteration.
\end{theorem}
\begin{proof}
    From Theorem \ref{thm...quadconv} we have
\begin{align}\label{eq...epsilonk}
	\varepsilon_{k+1}\leq  C\varepsilon_k^2\ \ \Rightarrow\ \ \varepsilon_k\leq (\varepsilon_0C)^{2^k}C^{-1}, \quad \text{for\ }k\in \mathbb{N}.
\end{align}
Consider the initial error $\tilde{\varepsilon} = \frac{1}{2^n}\varepsilon_0$ where $n=0,1,2,...$, and we replace $\varepsilon_0$ in \eqref{eq...epsilonk} to be $\tilde{\varepsilon}$. Then we have
\begin{align}
    \varepsilon_k\leq(\frac{1}{2^n}\varepsilon_0C)^{2^k}C^{-1}.
\end{align}
This implies that the Newton iteration converges if and only if $\frac{1}{2^n}\varepsilon_0C<1$, i.e. $n>\log_2(\varepsilon_0 C)$. 

Then an upper bound of the total number of iterations, $M_{iter}$, is an integer with $(\frac{1}{2^n}\varepsilon_0C)^{2^{M_{iter}}}C^{-1}<\varepsilon_{tol}$, i.e.
\begin{align}\label{eq...numiter}
    2^{M_{iter}}\log_2(\frac{1}{2^n}\varepsilon_0C)<\log_2(\varepsilon_{tol}C).
\end{align}
If $\varepsilon_{tol}C\geq1$, then the improved initial estimate $\tilde{\varepsilon}$ satisfies the tolerance. Thus in the following, we only consider the case where $\varepsilon_{tol}C<1$. From \eqref{eq...numiter} we know
\begin{align*}
    2^{M_{iter}}|\log_2(\frac{1}{2^n}\varepsilon_0C)|>|\log_2(\varepsilon_{tol}C)|
\end{align*}
since both $\log_2(\frac{1}{2^n}\varepsilon_0C)$ and $\log_2(\varepsilon_{tol}C)$ are negative. Hence
\begin{align*}
    M_{iter}>\log_2|\log_2(\varepsilon_{tol}C)|-\log_2|\log_2(\frac{1}{2^n}\varepsilon_0C)|
\end{align*}
is required such that the tolerance is satisfied. Remember that $M_{iter}$ should be nonnegative, so we choose
\begin{align*}
    M_{iter}=\max\{0,\ceil*{\log_2|\log_2(\varepsilon_{tol}C)|-\log_2|\log_2(\varepsilon_0C)-n|}+1\}.
\end{align*}
By denoting $C_1 = \log_2|\log_2(\varepsilon_{tol}C)|+1$ and $C_2=\log_2(\varepsilon_0C)$, we finally obtain that
\begin{align*}
    M_{iter}=\max\{0, \ceil*{C_1-\log_2|C_2-n|}\}
\end{align*}
holds for $n\in\mathbb{N}$ and $n>\log_2(\varepsilon_0 C)=C_2$.
\end{proof}

\begin{remark}
    We first provide a remark on \eqref{eq...epsilonk}. In \eqref{eq...epsilonk}, the factor $\varepsilon_0^{2^k}$ indicates that a small $\varepsilon_0$ can significantly reduce the error of the $k$-th step iteration $\varepsilon_k$. Therefore, the importance of initialisation in Newton iteration is demonstrated by this error estimate. It is worth noting that in general the exact Newton iteration count is not a monotonically decreasing function of $\varepsilon_0$. And this fact is not in conflict with our conclusion in Theorem \ref{thm...asym_iteration_count} because we only provide an upper bound of the iteration count. The exact iteration count should be smaller than this upper bound but not necessarily a monotonically decreasing function of $\varepsilon_0$.
\end{remark}

\begin{remark}
To enhance understanding of Theorem \ref{thm...quadconv}, here we provide some comments pertaining to the constant $C=\frac{1}{2} N^2 A_G B_G$ in \eqref{eq:C}. We consider solving the Allen--Cahn equation using the midpoint method as described in Section \ref{sec:numerical_examples}. Based on the definitions of $A_G$ and $B_G$, we find 
$$
C\propto N^2\tau\|u\|
$$
where $\tau$ is the time step in the midpoint method and $\|u\|$ is a certain norm of the solution $u$. Combined with the estimation of total iteration count $M_{iter}$ in \eqref{eq...numiter}, it can be concluded from this proportional relation that (i) larger time steps require a larger number of Newton iterates; (ii) more spatial discretisation points require a larger number of Newton iterates.
\end{remark}

\begin{remark}
To further demonstrate our theoretical result in Theorem \ref{thm...asym_iteration_count}, we investigate the change of iteration count when decreasing the $L^2$ error of the initial guess of Newton iteration in our numerical examples in Section \ref{sec:numerical_examples}. Figure \ref{fig:1d initial error} shows the result of the one-dimensional case. The blue line with cross-markers exhibits the true iteration count in the experiment when initial error $\tilde{\varepsilon} = \frac{1}{2^n}\varepsilon_0, n=0,1,2,...,17$. The red dashed line is the graph of $C_1-\log_2|C_2-n|$ as a function of $n$. 
Figure \ref{fig:2d initial error} shows the result of the two-dimensional case. The blue cross line is the iteration count in the experiment when initial error $\tilde{\varepsilon} = \frac{1}{2^n}\varepsilon_0, n=0,1,2,...,12$, while the red dashed line is still the theoretical result. 
From the comparison between the blue and red lines, we find the experiment results are perfectly consistent with our theoretical upper bound both in these two cases.
 \begin{figure}[htbp]
	\centering
	\subfigure[]{
		\begin{minipage}[b]{0.45\textwidth}
			\includegraphics[width=1\textwidth]{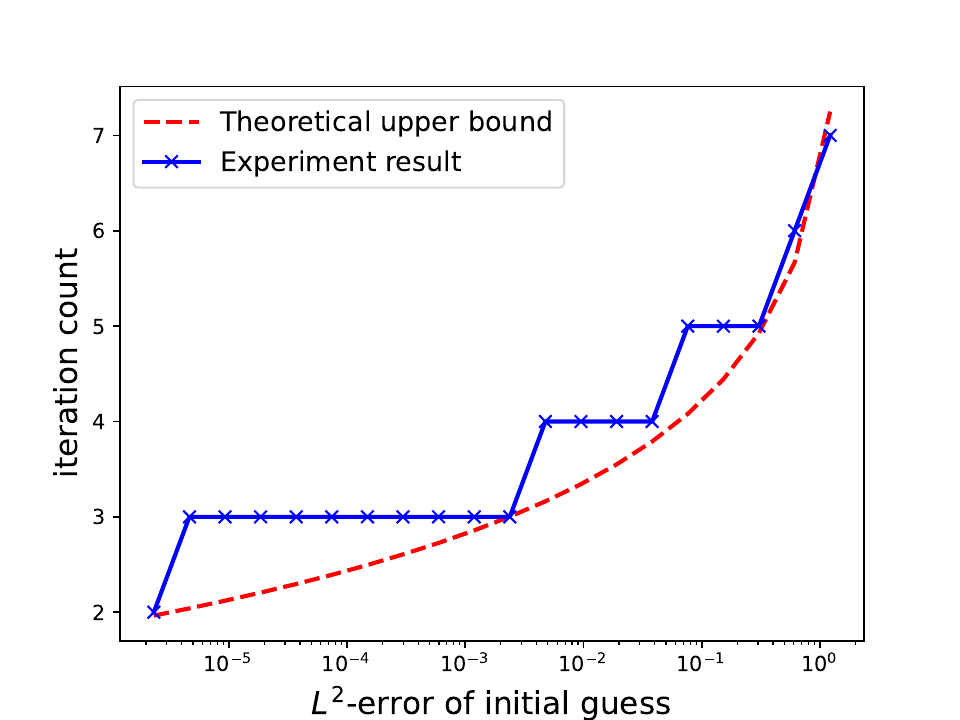}
		\end{minipage}
		\label{fig:1d initial error}
	}
    	\subfigure[]{
    		\begin{minipage}[b]{0.45\textwidth}
   		 	\includegraphics[width=1\textwidth]{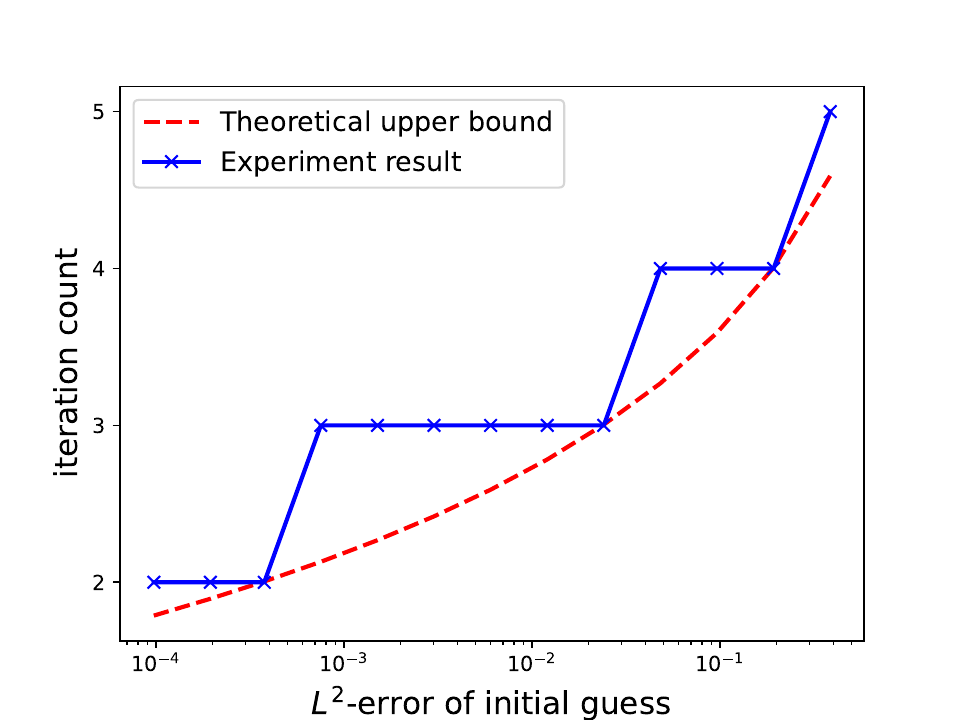}
    		\end{minipage}
		\label{fig:2d initial error}
    	}
	\caption{(a) Experimental and theoretical Newton iteration count in solving $1$-dimensional Allen--Cahn equation. The blue cross line is the experiment result, while the red dashed line is the theoretical estimation; (b) Experimental and theoretical Newton iteration count in solving $2$-dimensional Allen--Cahn equation. The blue cross line is the experiment result, while the red dashed line is the theoretical estimation.}
	\label{fig:initial error}
\end{figure}
\end{remark}

\subsection{Generalisation Error of implicit-scheme-informed Neural Networks}\label{sec:generalisation_error}

Having studied the cost reduction that can be achieved in Newton's method using  improved initialisations, we now turn our attention to understanding how good of an initialisation our neural network approximations (as described in Section~\ref{sec:hybrid_methods}) can provide. Our means of understanding the quality of our neural network initialisation will be to study the so-called `generalisation error' which we formally define in Sections~\ref{sec:fully-discrete_generalisation_error} \& \ref{sec:semi-discrete_generalisation_error}. Essentially this quantity provides a measure for how well we expect the initialisation to perform on unseen data, given knowledge of its performance during training. Before we can state the mathematical results more rigorously, we begin with introducing some useful concepts and notations: For a given Banach space $X$ with associated norm $\|\,\cdot\,\|_{X}$ and $z\in X$ we denote by $B_X(z,R)$ the ball of radius $R$ centred at $z$, i.e.
\begin{align*}
	B_{X}(z,R)=\{z\in X\vert \|z\|_{X}\leq R\}.
\end{align*}
We recall from Section \ref{sec:hybrid_methods}, that the exact solution of the implicit midpoint rule is denoted by $\Phi_\tau$ and then we denote the neural network approximation of $\Phi_{\tau}$ by $\Phi_{\tau}^*$. Recall that the definition of the midpoint rule can be expressed as $\Psi_\tau(u,\Phi_{\tau}(u))=\Phi_{\tau}(u)$. It turns out that a convenient way of studying the generalisation error for tensor product domains ($\Omega=\mathbb{T}^d=[-\pi,\pi]^d$) is via the framework of periodic Sobolev spaces and their appropriate subspaces for periodic, Dirichlet and Neumann boundary conditions respectively. In particular, for equations \eqref{eq...evol} with periodic boundary conditions we consider the Hilbert spaces
\begin{align*}
    H^{s}_{\mathrm{per}}(\Omega):=\{u\in L^2(\Omega)\,\vert\, \|u\|_{H^s_{\mathrm{per}}}<\infty\},
\end{align*}
where the periodic Sobolev norm is defined in the usual way. Given the Fourier coefficients
$$
\hat{u}_{\textbf{k}}=\frac{1}{(2\pi)^d}\int_{\Omega}u(x)\e^{-i\textbf{k}\cdot x} dx,
$$
we define
\begin{align*}
    \|u\|_{H^s_{\mathrm{per}}}^2:=\sum_{\textbf{k}\in\mathbb{Z}^d}\langle \textbf{k}\rangle^{2s}|\hat{u}_\textbf{k}|^2, \quad \text{where } \langle \textbf{k}\rangle:=\left\{ \begin{array}{cc}
     |\textbf{k}|\quad&\text{if } \textbf{k}\neq\textbf{0},\\
      1\quad&\text{if } \textbf{k}=\textbf{0}.\\
\end{array}
\right.
\end{align*}
In the following we use $H^s$ to replace $H^s_{\mathrm{per}}$ for notational simplicity. There are two central properties of these spaces that we exploit in what follows. The first lemma is the Sobolev embedding theorem which comes from Theorem 6.3 in \cite{1975}.
\begin{lemma}\label{lemma...sobemb}
    For any $\alpha>\beta\geq0$, the following embedding is compact:
    $$
    H^{\alpha}\rightarrow H^{\beta}.
    $$
    And we use the following notation to represent this compact embedding:
    $$
    H^{\alpha}\Subset H^{\beta}.
    $$
\end{lemma}
The second lemma is the well-known bilinear estimate.
\begin{lemma}\label{lemma...sobmulti}
    If $\alpha>d/2$, then for any $f,g\in H^\alpha(\mathbb{T}^d)$ the following inequality holds
    $$
    \|fg\|_{H^{\alpha}}\leq C_\alpha \|f\|_{H^\alpha}\|g\|_{H^\alpha}
    $$
    for some constant $C_\alpha>0$.
\end{lemma}

In case of Dirichlet boundary conditions let us assume without loss of generality that the domain is $\Omega=[0,\pi]^d$. Then we replace the aforementioned Sobolev spaces on $[-\pi,\pi]^d$ by subspaces arising from the condition
\begin{align*}
    \hat{u}_{\textbf{k}}=(-1)^{\sum_{l=1}^di_l}\hat{u}_{\textbf{j}},\,\, \text{whenever}\,\, \textbf{j}_l=(-1)^{i_l}\textbf{k}_{l},\,\, \text{for some } i_l\in\{0,1\}, l=1,\dots, d.
\end{align*}
This is equivalent to the subspaces spanned by Fourier sine series, i.e. (for $s>d/2$) those subspaces satisfying homogeneous Dirichlet conditions on $[0,\pi]^d$. For homogeneous Neumann conditions on $[0,\pi]^d$ we can proceed analogously by considering those subspaces of $H^s_{\mathrm{per}}([-\pi,\pi]^d)$ consisting of functions which satisfy
\begin{align*}
    \hat{u}_{\textbf{k}}=\hat{u}_{\textbf{j}},\,\, \text{whenever}\,\, \textbf{j}_l=(-1)^{i_l}\textbf{k}_{l},\,\, \text{for some } i_l\in\{0,1\}, l=1,\dots, d,
\end{align*}
i.e. subspaces spanned by Fourier cosine series. Clearly both of those choices of subspaces inherit the properties described in Lemmas \ref{lemma...sobemb} \& \ref{lemma...sobmulti} allowing us to state and prove analogous results as shown below in the periodic case. In the interest of brevity we will therefore focus on the periodic setting for the rest of this section. 

While we expect several of the statements in this section to generalise to a broader setting we choose to present our theoretical analysis for a specific equation and numerical integrator and spatial domains of dimension $d=1,2,3$. In particular, we consider the Allen--Cahn equation given as follows
\begin{equation}\label{eq...ac}
\left\{\begin{array}{l}
\partial_t u(\textbf{x},t) = \varepsilon_{AC}^2 \Delta u(\textbf{x},t) - f(u(\textbf{x},t)), t\in [0,T], \textbf{x}\in\Omega,\\
u(\textbf{x},0) = u_0(\textbf{x}), \textbf{x}\in \overline{\Omega},\\
\end{array}\right.
\end{equation}
where $\Omega=[-\pi,\pi]^d (d = 1,2,3)$ (or $\Omega=[0,\pi]^d$ with Dirichlet or homogeneous Neumann boundary conditions), $\varepsilon_{AC}>0$ is the interfacial width, and the nonlinear function $f(u) = u^3-u$. Our implicit method of choice for this section is the implicit midpoint scheme given by
\begin{equation}\label{scheme...midpoint}
u^{n+1}=\Psi_{\tau}(u^n,u^{n+1}) := u^n + \tau \mathcal{L}(\frac{1}{2}(u^n+u^{n+1}))
\end{equation}
where $\mathcal{L}(u) = \varepsilon_{AC}^2 \Delta u - u^3 + u$ and $\tau$ is the time step.
\begin{remark}\label{rmk...epsilon}
    To avoid misunderstanding, we clarify several symbols whose notations appear similar but have different meanings. Specifically, $\varepsilon_{AC}$ refers to the interfacial width, a parameter in Allen-Cahn equation \eqref{eq...ac}. $\varepsilon_k$ is the error of the k-th iteration in Newton's method as stated in Theorem \ref{thm...quadconv}. On the other hand, $\epsilon$ means a small positive number.
\end{remark}

The training of our neural network is performed over a set $\mathcal{T}$ of $N_{data}$ training data points which are, as indicated in Section \ref{sec:hybrid_methods}, of specific differentiability and uniformly bounded on the domain of interest, i.e. we assume for all $1\leq j\leq N_{data}$, $u_j,\Phi_{\tau}(u_j)\in H^\alpha(\Omega)$ for some $\alpha>d/2+2$.  In particular we sample training data from a bounded ball in $H^\alpha(\Omega)$ which, without loss of generality, we take to be the unit ball for the purposes of this analysis, i.e. $u_j\in B_{H^{\alpha}}(0,1)$. From the maximum bound principle (MBP) of the Allen--Cahn equation (see Section 
\ref{sec:numerical_examples} for details), it makes sense to choose uniformly bounded training data of the form $u_j\in B_{L^\infty}(0,1)$. In summary, we will prove an upper bound on the generalisation error for training data of the form
\begin{align}\label{eq...training_data}
	\mathcal{T}=\{u_j\}_{j=1}^{N_{data}}\subset B_{H^{\alpha}}(0,1)\cap  B_{L^\infty}(0,1),\quad \text{\ for\ }\alpha>d/2+2.
\end{align}

In the following we will provide two ways of studying the performance of the neural network initialisation on an unseen function. Firstly, we will consider the fully discrete case, taking into account the finite-dimensional nature of the input and output in our implementation of the neural network. Secondly, we will study the formal limit as the dimension of input and output of our neural network initialisation tend to infinity (i.e. considering a generic neural network with the scheme-informed loss \eqref{eq...loss}). We will see that even in this formal limit we can ensure that a finite number of training data points is sufficient for a close relation between the a-posteriori maximum training error and the maximum generalisation error over a suitably bounded subspace of functions.

In order to make a statement about the generalisation error we assume that the implicit midpoint rule for the Allen--Cahn equations applied to functions of the given regularity has a unique solution. This is similar to Assumption~\ref{assump...G} but now in the continuous case.
\begin{assumption}\label{asump...mainthm} {The implicit midpoint scheme always has a unique and bounded solution, i.e. given $u\in H^\alpha$, the solution $\Phi_{\tau}(u)$ uniquely exists and is bounded in $\|\cdot\|_{H^\alpha}$.}
\end{assumption}

\subsubsection{Fully discrete generalisation error analysis}\label{sec:fully-discrete_generalisation_error}
Having established the above basis, we begin by considering the fully discrete error analysis. This means that we look at our neural network initialisation as a map on finite-dimensional vector spaces (cf. our implementation in Section~\ref{sec:numerical_examples}):
\begin{align*}
\Phi_{\tau}^*:\mathbb{C}^{(2M+1)^d}&\mapsto\mathbb{C}^{(2M+1)^d}\\
(u^n(\mathbf{x}_{\mathbf{j}}))_{|\mathbf{j}|_{\infty}\leq M}&\mapsto (\Phi_{\tau}^*(\mathbf{u}^n)_{\mathbf{j}})_{|\mathbf{j}|_{\infty}\leq M},
\end{align*}
with $\Phi_{\tau}^*(\mathbf{u}^n)_{\mathbf{j}}\approx \Phi_{\tau}(u^n)(\mathbf{x}_{\mathbf{j}}),$  $\mathbf{x}_{\mathbf{j}}=2\pi\mathbf{j}/(2M+1)$ and $ \mathbf{j}\in\mathbb{Z}^{d}$, where $\Phi_{\tau}(u^n)$ is the midpoint rule \eqref{scheme...midpoint}. We then introduce the finite-dimensional subspaces of trigonometric polynomials, $S_M$,
\begin{align*}
    S_M:=\left\{f(\mathbf{x})=\sum_{\substack{\mathbf{k}\in\mathbb{Z}^d\\|\mathbf{k}|_\infty\leq M}}c_{\mathbf{k}}e^{i\mathbf{x}\cdot \mathbf{k}}\,\Big\vert\, c_{\mathbf{k}}\in\mathbb{C}, \text{for\ all\ }|\mathbf{k}|_{\infty}\leq M\right\},
\end{align*}
where for a multi-index $\mathbf{k}\in\mathbb{Z}^d$ we denote $|\mathbf{k}|_{\infty}=\max_{1\leq j\leq d}|k_j|$. This allows us to equivalently look at the action of $\Phi^*_{\tau}$ on $S_M$
 using the following encoding and recovery operators which allow for an identification of the finite dimensional space $S_M$ with $\mathbb{C}^{(2M+1)^d}$ and take the form:\\
\begin{minipage}{0.5\textwidth}
\begin{align*}
    \mathbf{{R}}:\quad \mathbb{C}^{(2M+1)^d}&\rightarrow S_{M},\\
                (u(\mathbf{x}_{\mathbf{j}}))_{|\mathbf{j}|_{\infty}\leq M}&\mapsto u,
\end{align*}
\end{minipage}%
\begin{minipage}{0.5\textwidth}\begin{align*}\mathbf{{E}}: S_{M}&\rightarrow \mathbb{C}^{(2M+1)^d},\\
                u&\mapsto (u(\mathbf{x}_{\mathbf{j}}))_{|\mathbf{j}|_{\infty}\leq M},
\end{align*}
\end{minipage}\\\ \\
where $\mathbf{j}\in\mathbb{Z}^{d}, |\mathbf{j}|_\infty\leq M$. These operators are simply discrete Fourier transforms and are isometries between $\mathbb{C}^{(2M+1)^d}$ with the Euclidean norm and $S_M$ with the $L^2$-norm. {We note that such a fully discrete framework is the commonly used way of analysing the performance of oprator-learning frameworks (for example in the approximation error of PCA autoencoders \cite{finite_dim_appro_error}, approximation and generalization errors for DeepONets \cite{DeepOnet_error}, and the generalisation error for FNO \cite{FNO_error}}).
\begin{remark}
    As noted above in the case of Dirichlet, or Neumann boundary conditions we can perform a similar type of analysis and equivalence of finite-dimensional spaces by looking instead at Fourier sine or cosine series respectively.
\end{remark}

After training the neural network with the loss \eqref{eq...loss} we have an automatic control on the average training error $\bar{\mathcal{E}}_{T}$, which in terms of the training dataset $\mathcal{T}$ as given in \eqref{eq...training_data} is:
\begin{align}\label{eq...avg_trainerror}
	\bar{\mathcal{E}}_{T}&=\frac{1}{N_{{data}}}\sum_{j=1}^{N_{data}}\|\Psi_\tau(\mathbf{{R}}(\mathbf{c}_j),\mathbf{{R}}\circ\Phi^*_{\tau}(\mathbf{c}_j))-\mathbf{{R}}\circ\Phi^*_{\tau}(\mathbf{c}_j)\|_{L^2}^2,
\end{align}
In the following, we will estimate the error due to our implicit-scheme-informed learning strategy in approximating $\Phi_\tau$ for unseen data in the relevant set. In particular, we consider the average error over $\mathbf{E}(B_{H^\alpha}(0,1)\cap  S_M)$, also referred to as the average generalisation error given by:
\begin{align*}
    \bar{\mathcal{E}}_{G}&:={\int_{\mathbf{{E}}(S_{M}\cap  B_{H^\alpha}(0,1))}\hspace{-2.675cm}\mathbf{{-}}\hspace{2.675cm}\|\Phi_{\tau}\circ\mathbf{{R}}(\mathbf{z})-\mathbf{{R}}\circ\Phi^*_{\tau}(\mathbf{z})\|_{L^2}^2 \,d\mathbf{z}}\\
    &=\frac{\int_{\mathbf{{E}}(S_{M}\cap  B_{H^\alpha}(0,1))}\|\Phi_{\tau}\circ\mathbf{{R}}(\mathbf{z})-\mathbf{{R}}\circ\Phi^*_{\tau}(\mathbf{z})\|_{L^2}^2 \,d\mathbf{z}}{\int_{\mathbf{{E}}(S_{M}\cap  B_{H^\alpha}(0,1))} \,d\mathbf{z}}.
\end{align*}

\begin{theorem}\label{thm:averaged_generalisation_error}  Let $d\in \{1,2,3\}$, fix $M\in \mathbb{N}$, $\alpha> 2$ and suppose that Assumption \ref{asump...mainthm} is satisfied. Then there exist a constant $C=C(M,\alpha,d)>0$ such that the following holds. Given $\epsilon>0$ there is a $N^{(M,d)}_\epsilon\in\mathbb{N}$ such that for all $N_{data}>N^{(M,d)}_\epsilon$ we can find a finite sequence of training data $\{u_j\}_{j=1}^{N_{data}}$ such that for all $\tau<\tau_0$ we have the estimate
    \begin{align}\label{eqn:estimate_averaged_generalisation_error}
        \bar{\mathcal{E}}_{G}\leq C(\epsilon+\bar{\mathcal{E}}_{T}).
    \end{align}
    Furthermore, $N_\epsilon^{(M,d)}$ scales no worse than
    \begin{align*}
        N_\epsilon^{(M,d)}\lesssim  \left\lceil\left(\frac{2}{\epsilon}\right)^{2^d(2M+1)^d}\right\rceil.
    \end{align*}
\end{theorem}
\begin{proof}
The proof of Theorem \ref{thm:averaged_generalisation_error} is given in Section \ref{sec:proof_of_main_thm}.
\end{proof}
Theorem~\ref{thm:averaged_generalisation_error} provides us with a guarantee that, for any fixed spatial resolution $M$, we can find a sufficient number of training data points such that (up to a constant factor depending on the Lipschitz constant of the neural network) the generalisation error is arbitrarily close to the training error.
    
\subsubsection{Semi-discrete generalisation error analysis}\label{sec:semi-discrete_generalisation_error}
We note that in the above, fully-discrete case, the number of training data required to achieve $\epsilon$-closeness in the averaged generalisation error has to grow exponentially in $M$, as does the constant $C$ in \eqref{eqn:estimate_averaged_generalisation_error}. It turns out that by looking at the maximum generalisation error we can extend this result to an estimate on the full infinite-dimensional case (i.e. the formal limit $M\rightarrow\infty$) while ensuring the corresponding constants in the estimates remain finite. In this formal limit we will regard our neural network as a map on function spaces:
\begin{align*}
    \Phi^*_{\tau}:L^2\rightarrow L^2,
\end{align*}
where $\Phi^*_{\tau}(u^n)$ is still the approximation of midpoint rule $\Phi_{\tau}(u^{n})$ defined in \eqref{scheme...midpoint}. This point-of-view for $\Phi^*_{\tau}$ can be interpreted in two ways: (i) as an operator which is the formal limit of a sequence of fully discrete neural networks taking finite-dimensional inputs and outputs as described in the previous section; (ii) a study of the generalisation error of an arbitrary operator learning architecture (which takes in functions and outputs functions) when used with the scheme-informed loss \eqref{eq...loss}.

While the averaged training error \eqref{eq...avg_trainerror} is used as the optimisation objective in the neural network training, once the training is completed we have access to the maximum training error $\mathcal{E}_{T}^{\infty}$, which for mathematical convenience we will express in this section in terms of the maximum over the training dataset $\mathcal{T}$ as given in \eqref{eq...training_data}:
\begin{align*}
	\mathcal{E}_{T}^{\infty}&=\max_{1\leq j\leq N_{data}}\|\Psi_\tau(u_j,\Phi^*_{\tau}(u_j))-\Phi^*_{\tau}(u_j)\|_{L^2}.
\end{align*}
In the following, we will estimate the error due to our implicit-scheme-informed learning strategy in approximating $\Phi_\tau$ for unseen data in the relevant set. In particular, we consider the total error over $B_{H^{\alpha}}(0,1)\cap  B_{L^\infty}(0,1)$, also referred to as the maximum generalisation error given by:
\begin{align*}
    \mathcal{E}_{G}^{\infty}=\sup_{\substack{u\in B_{H^\alpha}(0,1)}\cap B_{L^\infty}(0,1)}\|\Phi_{\tau}(u)-\Phi^*_{\tau}(u)\|_{L^2}.
\end{align*}

The following theorem states that for a sufficiently large training dataset, the generalisation error $\mathcal{E}_{G}^{\infty}$ is arbitrarily close to the training error $\mathcal{E}_{T}^{\infty}$. Notably the constants in this estimate are finite, even though we work on infinite-dimensional spaces, thus avoiding the exponential growth (i.e. the curse of dimensionality) encountered in the constants of Theorem~\ref{thm:averaged_generalisation_error}.

\begin{theorem}\label{thm...genlerror} Suppose $d\in\{1,2,3\}, \alpha> 2$ and that Assumption \ref{asump...mainthm} is satisfied. Then there exist constants $\tilde{C}_1$, $\tilde{C}_2$ and $\tau_0$ such that the following holds. Given $\epsilon>0$ there is a $N_\epsilon\in\mathbb{N}$ such that for all $N_{data}>N_\epsilon$ we can find a finite sequence of training data $\{u_j\}_{j=1}^{N_{data}}$ such that for all $\tau<\tau_0$ we have the estimate
    \begin{align*}
        \mathcal{E}_{G}^{\infty}\leq \tilde{C}_1\mathcal{E}_{T}^{\infty}+\tilde{C}_2\epsilon.
    \end{align*}
    Furthermore, $N_\epsilon$ scales no worse than
    \begin{align*}
        N_{\epsilon}\lesssim \left\lceil\left(\frac{4}{\epsilon}\right)^{2^d\left(2(\frac{\epsilon}{2})^{\frac{1}{2}-\alpha}+1\right)^d}\right\rceil.
    \end{align*}
\end{theorem}
\begin{proof} The proof of Theorem \ref{thm...genlerror} is given in Section \ref{sec:proof_of_main_thm}.
\end{proof}
\begin{remark}
    For a better understanding of Theorem \ref{thm:averaged_generalisation_error} and \ref{thm...genlerror}, we further provide some remarks on them.
    \begin{itemize}
        \item The loss function \eqref{eq...loss} does not directly control the maximal training error when $N_{data}\rightarrow \infty$. While we recall that the average training error $\bar{\mathcal{E}}_T$ can be automatically controlled by the loss function \eqref{eq...loss} for any $N_{data}$.
        \item In practice, however, the maximum training error $\mathcal{E}_T^\infty$ can be easily calculated after the neural network is trained. Thus, we highlight that there is no a-priori control on $\mathcal{E}_T^\infty$ but we can evaluate it a-posteriori.
        \item We note that there are a significant number of references using advanced tools in statistical learning theory to bound the generalisation error, such as Rademacher complexity and VC dimension. However, such statistical learning framework requires the appropriate choice of the hypothesis space to balance the approximation error and the estimation error (see \cite{CMS_weinan} for detailed explanations), which is a delicate problem and beyond the scope of this paper. In contrast to statistical learning analysis, our estimates presented in Theorems \ref{thm:averaged_generalisation_error} \& \ref{thm...genlerror} are completely deterministic and provide a link between the generalisability of such models directly to the function spaces relevant to the analysis of classical numerical methods. This is possible mainly because we focus on a specific PDE at hand, and provides a justification of our scheme-informed (unsupervised) loss for the initialisation task in Newton's method: we do not need to generate training data-pairs (i.e. pairs $(u,\Phi_{\tau}(u))$) and still we can claim that (under appropriate standard assumptions on the Lipschitz continuity of the neural network) the method generalises well to unseen data if the training data are appropriately selected. 
        The independence of our analysis from the specific neural network architecture at hand also allows for the potential expansion of the proposed hybrid solver to other equations with enough regularity of solutions, implicit schemes with stability and network structures that are Lipschitz continuous.
        \item While deterministic and probabilistic estimates share the common goal of bounding the difference between training and generalization errors, unlike high-probability statements, which assume that the training data are drawn i.i.d. from a probability distribution, our approach focuses on the existence of such a sequence
of training data. Replacing the deterministic estimate with a high-probability statement via concentration bounds may be feasible in the fully-discrete setting but is nontrivial due to the implicit-scheme-based residual in the training error. Deriving appropriate concentration bounds remains an interesting direction for future research.
    \end{itemize}
\end{remark}

 \subsubsection{Minimal number of training data required for maximum bound on the generalisation error}
Finally, we will investigate how much training data are required in order to achieve the accuracy $\epsilon>0$. 
 From the the second part of the proof of Theorem \ref{thm...genlerror} (see Section \ref{sec:proof_of_main_thm} for details), we find the amount of training data required is equivalent to the size of an $\epsilon$-net covering $B_{H^\alpha}(0,1)\cap B_{L^\infty}(0,1)$ in $H^2$. The following proposition provides a quantifiable estimate on the size of such an $\epsilon$-net in $H^{\gamma}$.
\begin{proposition}[Amount of training data required to achieve $\epsilon$-accuracy.]\label{prop...num_traindata}
Fix $\alpha>\gamma>d/2.$ Given $\epsilon>0$ there is an $\epsilon$-Net of size at most
\begin{align*}
    N_{\epsilon}^{(M,d)}=\left\lceil\left(\frac{2}{\epsilon}\right)^{2^d(2M+1)^d}\right\rceil
\end{align*}
covering $S_M\cap B_{H^\alpha}(0,1)\cap B_{L^\infty}(0,1)$ in $H^{\gamma}$. Moreover, there is an $\epsilon$-Net of size at most
\begin{align*}
    N_\epsilon=\left\lceil\left(\frac{4}{\epsilon}\right)^{2^d\left(2(\frac{\epsilon}{2})^{\frac{1}{\gamma}-\alpha}+1\right)^d}\right\rceil
\end{align*}
covering $B_{H^\alpha}(0,1)\cap B_{L^\infty}(0,1)$ in $H^{\gamma}$.
\end{proposition}
\begin{proof} The proof can be found in Appendix \ref{appendix:proof_prop}.
\end{proof}

\subsection{Proof of the main theorems}\label{sec:proof_of_main_thm}
In this section, we show the detailed proofs of the two main results, Theorems \ref{thm:averaged_generalisation_error} and \ref{thm...genlerror}. We begin by introducing two auxiliary estimates which will be used in the proofs. In the following two lemmas the Lipschitz function $G$ should be ultimately thought of as the neural network $\Phi^*_\tau$ as will be apparent in the proof of main theorems.

Our first step in both types of analysis is to relate the error in an unseen datapoint $u$ to the residual error in the scheme-informed loss (cf. \eqref{eq...loss}).
\begin{lemma}\label{lem:loss_bounds_approximation_error}
    Suppose $\alpha>\gamma\geq 0$ and $\alpha>d/2$ and that Assumption~\ref{asump...mainthm} is satisfied, and let $G:H^{\gamma}\rightarrow H^{\gamma}$ be a Lipschitz continuous map. Then there is a $\tau_0>0$ and a constant $C>0$ such that for all $u\in B_{H^{\alpha}}(0,1)\cap  B_{L^\infty}(0,1)$ and all $0<\tau\leq \tau_0$ the following estimate holds:
    \begin{align*}
         \|\Phi_{\tau}(u)-G(u)\|_{H^{\gamma}}\leq C\|\mathcal{R}(G(u))\|_{H^\gamma},
\end{align*}
where $\mathcal{R}(G(u)) = \mathcal{R}(G(u);u)=\Psi_\tau(u,G(u))-G(u)$.
\end{lemma}
\begin{proof}
    The proof is provided in Appendix \ref{appendix:two_lemma}.
\end{proof}

The next result provides a useful estimate of the residual (scheme-informed error) in terms of the distance between datapoints, which will be used in the proofs of the main theorems in the following sections.

\begin{lemma}\label{lem:estimate_onresidual}
Suppose $\alpha> \beta+2>d/2$ and that Assumption \ref{asump...mainthm} is satisfied. Let $G:H^{\beta}\rightarrow H^\beta$ be a Lipschitz continuous map, then there are constants $C,\tau_0>0$ such that for all $0\leq\tau\leq \tau_0$ and any $v,w\in B_{H^\alpha}(0,1)$ we have
\begin{align}
    \|\mathcal{R}(G(v))-\mathcal{R}(G(w))\|_{H^\beta}\leq C\|v-w\|_{H^{\beta+2}}.
\end{align}
\end{lemma}
\begin{proof}
    The proof is provided in Appendix \ref{appendix:two_lemma}.
\end{proof}

After having the above two lemmas, we first provide the proof of Theorem \ref{thm:averaged_generalisation_error} as follows.
\begin{proof}(Proof of Theorem \ref{thm:averaged_generalisation_error})
To begin with, we observe that, similar to the steps in \cite{DBLP:journals/corr/SzegedyZSBEGF13}, it is easy to show our neural network is Lipschitz continuous. Thus $\mathbf{{R}}\circ\Phi_{\tau}^*\circ\mathbf{{E}}: L^2\rightarrow L^2$ is Lipschitz continuous and we have by Lemma~\ref{lem:loss_bounds_approximation_error} (since $1\leq d\leq 3$), for some constants $C_1>0$ and $\tau_0>0$ which depend on the Lipschitz constant of the trained neural network,
        \begin{align*}
            \left(\int_{\mathbf{{E}}(S_{M}\cap  B_{H^\alpha}(0,1))} \,\hspace{-.5cm}d\mathbf{z}\right)\bar{\mathcal{E}}_{G}&\leq C_1\int_{\mathbf{{E}}(S_{M}\cap  B_{H^\alpha}(0,1))}\|\Psi_\tau(u,\mathbf{{R}}\circ\Phi^*_{\tau}(\mathbf{z}))-\mathbf{{R}}
            \circ\Phi^*_{\tau}(\mathbf{z})\|_{L^2}^2 d\mathbf{z}.
        \end{align*}
Now we note that by Proposition~\ref{prop...num_traindata} we can find $\mathbf{c}_j\in \mathbf{{E}}(S_{M}\cap  B_{H^\alpha}(0,1)), j=1,\dots, N_{\epsilon}^{(M,d)}$ with 
\begin{align*}
    N_{\epsilon}^{(M,d)}\leq \left\lceil\left(\frac{2}{\epsilon}\right)^{2^d(2M+1)^d}\right\rceil,
\end{align*}
such that
\begin{align*}
    \mathbf{{E}}(S_{M}\cap  B_{H^\alpha}(0,1))\subset\bigcup_{j=1}^{N_{\epsilon}^{(M,d)}}\mathbf{{E}}(S_M\cap B_{H^{2}}(\mathbf{{R}}(\mathbf{c}_j),\epsilon)).
\end{align*}
        Thus, we have
        \begin{align*}
            \left(\int_{\mathbf{{E}}(S_{M}\cap  B_{H^\alpha}(0,1))} \,\hspace{-.5cm}d\mathbf{z}\right)\bar{E}_{G}&\leq \sum_{j=1}^{N_{\epsilon}^{(M,d)}}\int_{\mathbf{{E}}(S_M\cap B_{H^{2}}(\mathbf{{R}}(\mathbf{c}_j),\epsilon))}\|\Psi_\tau(u,\mathbf{{R}}\circ\Phi^*_{\tau}(\mathbf{z}))-\mathbf{{R}}
            \circ\Phi^*_{\tau}(\mathbf{z})\|_{L^2}^2 d\mathbf{z}.
        \end{align*}
Moreover, by Lemma~\ref{lem:estimate_onresidual} (noting that $1\leq d\leq 3$),
\begin{align*}
    &\int_{\mathbf{{E}}(S_M\cap B_{L^2}(\mathbf{{R}}(\mathbf{c}_j),\epsilon))}\|\Psi_\tau(u,\mathbf{{R}}\circ\Phi^*_{\tau}(\mathbf{z}))-\mathbf{{R}}
            \circ\Phi^*_{\tau}(\mathbf{z})\|_{L^2}^2 d\mathbf{z}\\
            &=\int_{\mathbf{{E}}(S_M\cap B_{H^{2}}(\mathbf{{R}}(\mathbf{c}_j),\epsilon))}\|\mathcal{R}(\mathbf{{R}}\circ\Phi^*_{\tau}\circ\mathbf{{E}})\circ \mathbf{{R}}(\mathbf{z})\|_{L^2}^2d\mathbf{z}\\
            &\leq C_{NN}\int_{\mathbf{{E}}(S_M\cap B_{H^{2}}(\mathbf{{R}}(\mathbf{c}_j),\epsilon))}\|\mathbf{{R}}(\mathbf{z})-\mathbf{{R}}(\mathbf{z}_j)\|_{H^{2}}^2d\mathbf{z}\\
            &\quad\quad\quad\quad+\|\mathcal{R}(\mathbf{{R}}\circ\Phi^*_{\tau}\circ\mathbf{{E}})\circ \mathbf{{R}}(\mathbf{z}_j)\|_{L^2}^2\int_{\mathbf{{E}}(S_M\cap B_{H^{2}}(\mathbf{{R}}(\mathbf{c}_j),\epsilon))}d\mathbf{z}\\
            &\leq \left(C_{NN}\epsilon+\|\mathcal{R}(\mathbf{{R}}\circ\Phi^*_{\tau}\circ\mathbf{{E}})\circ \mathbf{{R}}(\mathbf{z}_j)\|_{L^2}^2\right)\int_{\mathbf{{E}}(S_M\cap B_{H^{2}}(\mathbf{{R}}(\mathbf{c}_j),\epsilon))}d\mathbf{z}
\end{align*}
where $C_{NN}>0$ depends on the Lipschitz constant of $\Phi_{\tau}^{*}$. Now we note that $\dim S_{M}=2^d(2M+1)^d$ and thus, since the norm $\|\,\cdot\,\|_{H^2}$ is linear with respect to multiplication by a positive scalar and $\mathbf{E}$ is a linear map, we have, for any $\mathbf{c}_j\in\mathbb{C}^{(2M+1)^d}$,
\begin{align*}
    \frac{\int_{\mathbf{{E}}(S_{M}\cap  B_{H^{2}}(\mathbf{R}(\mathbf{c}_j),\epsilon))} d\mathbf{z}}{\int_{\mathbf{{E}}(S_{M}\cap  B_{H^{2}}(0,1))} d\mathbf{z}}=\epsilon^{2^d(2M+1)^d}.
\end{align*}
Thus, in summary, we find the desired estimate
\begin{align*}
    \bar{E}_{G}&\leq C\left(\epsilon +\bar{E}_{T}\right),
\end{align*}
where
\begin{align*}
    C\leq C_{NN}\frac{\int_{\mathbf{{E}}(S_{M}\cap  B_{H^{2}}(0,1))} d\mathbf{z}}{\int_{\mathbf{{E}}(S_{M}\cap  B_{H^{\alpha}}(0,1))} d\mathbf{z}}\left(2^{2^d(2M+1)^d}+1\right),
\end{align*}
and $C_{NN}>0$ depends on the Lipschitz constant of $\Phi_{\tau}^{*}$. In the above we used that
\begin{align*}
    \epsilon^{2^d(2M+1)^d}\left\lceil\left(\frac{2}{\epsilon}\right)^{2^d(2M+1)^d}\right\rceil\leq 2^{2^d(2M+1)^d}+1.
\end{align*}
\end{proof}

We then prove the Theorem \ref{thm...genlerror}:
\begin{proof} The proof of Theorem \ref{thm...genlerror} can be divided into two parts.

\noindent
\textbf{Step 1:} We first estimate $\mathcal{E}_{G}^{\infty}$ in terms of the residual based error 
\begin{align*}
    \sup_{u\in B_{H^\alpha}(0,1)\cap B_{L^\infty}(0,1)}\|\mathcal{R}(\Phi_{\tau}^*(u))\|_{L^2}.
\end{align*}
From Lemma \ref{lem:loss_bounds_approximation_error}, since the neural network $\Phi_\tau^*$ is Lipschitz continuous (cf. \cite{DBLP:journals/corr/SzegedyZSBEGF13}), there exist $\tau_1>0$ and $\tilde{C}_1>0$ such that for all $u\in B_{H^\alpha}(0,1)\cap B_{L^\infty}(0,1)$ and $0<\tau<\tau_1$ we have by taking $\gamma=0$ in Lemma \ref{lem:loss_bounds_approximation_error}
\begin{align*}
    \|\Phi_\tau(u)-\Phi_\tau^*(u)\|_{L^2}\leq \tilde{C}_1\|\Phi_\tau^*(u)\|_{L^2}.
\end{align*}
Therefore after taking supremum we have the following estimate of generalisation error $\mathcal{E}_{G}^{\infty}$ in terms of the residual based error, valid for any $0<\tau<\tau_1$:
\begin{align}\label{eq...estimate_step1}
    \mathcal{E}_{G}^{\infty}\leq \tilde{C}_1\sup_{u\in B_{H^\alpha}(0,1)\cap B_{L^\infty}(0,1)}\|\mathcal{R}(\Phi_{\tau}^*(u))\|_{L^2}.
\end{align}

\noindent
\textbf{Step 2:} Then we estimate 
$$
\sup_{u\in B_{H^\alpha}(0,1)\cap B_{L^\infty}(0,1)}\|\mathcal{R}(\Phi_{\tau}^*(u))\|_{L^2}
$$
in terms of the training error $\mathcal{E}_{T}^{\infty}$.

\noindent
Since $B_{H^\alpha}(0,1)\cap B_{L^\infty}(0,1)$ is bounded in $H^\alpha$, combined with the conclusion in Lemma \ref{lemma...sobemb} we know $B_{H^\alpha}(0,1)\cap B_{L^\infty}(0,1)$ is totally bounded in $H^{2}$. Thus for any given $\epsilon>0$ there is a finite $\epsilon$-net w.r.t. $\|\cdot\|_{H^{2}}$ covering $B_{H^\alpha}(0,1)\cap B_{L^\infty}(0,1)$. i.e. there is an $N_\epsilon\in\mathbb{N}$ and $u_j\in B_{H^\alpha}(0,1)\cap B_{L^\infty}(0,1)$ for $1\leq j\leq N_{\epsilon}$ such that
\begin{align*}
    B_{H^\alpha}(0,1)\cap B_{L^\infty}(0,1)\subset\bigcup_{j=1}^{N_{\epsilon}}B_{H^2}(u_j,\epsilon).
\end{align*}
An upper bound on the size of $N_\epsilon$ is given in Proposition \ref{prop...num_traindata} by taking $\gamma=2$. For simplicity of notation, we shall write $B_\cap=B_{H^\alpha}(0,1)\cap B_{L^\infty}(0,1)$ in what follows. We now have 
\begin{align*}
    \sup_{u\in B_\cap}\|\mathcal{R}(\Phi_{\tau}^*(u))\|_{L^2}=\max_{1\leq j\leq N_{\epsilon}}\sup_{v\in B_{H^2}(u_j,\epsilon)\cap B_\cap}\|\mathcal{R}(\Phi_{\tau}^*(v))\|_{L^2}.
\end{align*}
Now we have for any $v\in B_{H^2}(u_j,\epsilon)$
\begin{align*}
    \|\mathcal{R}(\Phi_{\tau}^*(v))\|_{L^2}\leq\|\mathcal{R}(\Phi_{\tau}^*(u_j))\|_{L^2}+\|\mathcal{R}(\Phi_{\tau}^*(v))-\mathcal{R}(\Phi_{\tau}^*(u_j))\|_{L^2}.
\end{align*}
From Lemma~\ref{lem:estimate_onresidual} and the Lipschitz continuity of $\Phi_\tau^*$, there are constants $C_2$, $\tau_2$ such that for all $0<\tau\leq\tau_2$ we have by taking $\beta=0$ in Lemma~\ref{lem:estimate_onresidual}
\begin{align*}
    \|\mathcal{R}(\Phi_{\tau}^*(v))-\mathcal{R}(\Phi_{\tau}^*(u_j))\|_{L^2}\leq C_2\|v-u_j\|_{H^2}\leq C_2 \epsilon,
\end{align*}
which follows that
\begin{align*}
     \|\mathcal{R}(\Phi_{\tau}^*(v))\|_{L^2}\leq\|\mathcal{R}(\Phi_{\tau}^*(u_j))\|_{L^2}+C_2 \epsilon.
\end{align*}
It follows immediately by taking supremum over $B_\cap$ that
\begin{align}\nonumber
    \sup_{u\in B_\cap}\|\mathcal{R}(\Phi_{\tau}^*(u))\|_{L^2}&=\max_{1\leq j\leq N_{\epsilon}}\sup_{v\in B_{H^2}(u_j,\epsilon)\cap B_\cap}\|\mathcal{R}(\Phi_{\tau}^*(v))\|_{L^2}\\\nonumber
    &\leq \max_{1\leq j\leq N_{\epsilon}}\|\mathcal{R}(\Phi_{\tau}^*(u_j))\|_{L^2}+C_2 \epsilon\\\label{eq...estimate_step2}
    &=\mathcal{E}_T^\infty+C_2 \epsilon.
\end{align}
This completes Step 2. 
Finally, combining \eqref{eq...estimate_step1} and \eqref{eq...estimate_step2} leads to the following estimate, valid for $0<\tau\leq\tau_0=\min\{\tau_1, \tau_2\}$,
\begin{align*}
    \mathcal{E}_G&\leq \tilde{C}_1\sup_{u\in B_\cap}\|\mathcal{R}(\Phi_{\tau}^*(u))\|_{L^2}\\
    &\leq \tilde{C}_1\mathcal{E}_T^\infty+\tilde{C}_2\epsilon
\end{align*}
where $\tilde{C}_2=\tilde{C}_1C_2$. This concludes the proof of Theorem \ref{thm...genlerror}.
\end{proof}

\section{Numerical examples and evaluation}\label{sec:numerical_examples}
Following the theoretical results in previous sections we now provide several numerical experiments to demonstrate the speed-up that can be achieved using our neural hybrid solver to initialise Newton's method. Aligned with the theoretical analysis in Section \ref{sec:scheme_informed_learning}, we consider solving the Allen--Cahn equation \eqref{eq...ac}, a reaction--diffusion equation that describes the motion of a curved antiphase boundary \cite{ALLEN19791085}, equipped with homogeneous Neumann boundary condition. It is known that equation \eqref{eq...ac} is the $L^2$ gradient flow with respect to the following energy functional
\begin{equation}\label{eq...energy}
E(u) = \int_{\Omega} (\frac{\varepsilon_{AC}^2}{2}|\nabla u|^2 + F(u)) \,d\textbf{x},
\end{equation}
where $F(u)=\frac{1}{4}(u^2-1)^2$ such that $F'(u)=f(u)$. Observe that
$$
\frac{d}{dt}E(u(t)) = -\int_{\Omega}|\partial_t u(\textbf{x},t)|^2 \,d\textbf{x}\leq 0.
$$
Hence, the energy \eqref{eq...energy} dissipates along the evolution of \eqref{eq...ac}.

Another important property of \eqref{eq...ac} is that its solution satisfies the maximum bound principle (MBP). Namely if 
$$
\max_{\textbf{x}\in\overline{\Omega}}|u_0(\textbf{x})|\leq 1,
$$
then the following also holds
$$
\max_{\textbf{x}\in\overline{\Omega}}|u(\textbf{x},t)|\leq 1 \quad\text{for all}\quad t>0.
$$
Energy dissipation and MBP are two vital physical properties of \eqref{eq...ac}. Thus, when constructing numerical methods, we hope to preserve them in order to achieve better long-term behaviour.

Numerous numerical methods are exploited in the literature for solving semilinear parabolic equations like \eqref{eq...ac}. \cite{gomez2011} proposed a provably unconditionally stable mixed variational method for the class of phase-field models. However, only the nonlinear stability restricted to the free energy can be guaranteed in this method. \cite{im-exp} explored an MBP-preserving method for the Allen--Cahn equation which uses implicit-explicit scheme in time and central finite difference in space. Exponential-type integrators have been found to have some advantages in solving equations containing strong stiff linear terms. An example of the exponential-type integrators is the exponential time differencing (ETD) schemes studied in \cite{doi:10.1137/19M1243750}. In ETD methods, which are based on Duhamel's forumla, the linear part is evaluated exactly leading to both stability and accuracy. As we noted in Section \ref{sec:intro}, fully implicit schemes, such as implicit midpoint methods, often exhibit better stability and accuracy with even relatively large time steps. However, implicit methods are more challenging to employ, due to their high computational cost.

Throughout this section we focus on solving \eqref{eq...ac} by the implicit midpoint method \eqref{scheme...midpoint}, the most common fully implicit method, in time and using a central finite difference discretisation in space. However, we note that our proposed method provides a general framework of hybrid Newton solvers for implict numerical methods that is not limited to the present discretisation. As a baseline for judging whether our proposed neural hybrid solver can accelerate the whole time-stepping algorithm and achieve a better performance we compare also against explicit methods, such as the ETD method. Before discussing the numerical results let us provide a few additional details on our implementation. First of all, we write down the central finite difference scheme to approximate the spatial derivatives as shown in \cite{im-exp}. Denote by $D_h$ the discrete matrix of the Laplacian $\Delta$ equipped with homogeneous Neumann boundary condition in $1$-dimensional case.
$$
D_h = \Lambda_h := \frac{1}{h^2} \begin{bmatrix}
-1 & 1 &  &  &  \\
1 & -2 & 1 &  &  \\
 & \ddots & \ddots & \ddots &  \\
 &  & 1 & -2 & 1 \\
 &  &  & 1 & -1 
\end{bmatrix} _{N\times N}
$$
where $h$ is the spatial increment. Let $\otimes$ denote the Kronecker product and $I$ is the identity matrix, then the discrete matrix of Laplacian in the $2$-dimensional case can be written as
$$
D_h = I\otimes\Lambda_h+\Lambda_h\otimes I.
$$

The implicit midpoint method has been given in \eqref{scheme...midpoint}. The time step $\tau$ here is also referred to as $\tau_{midpoint}$ in the following subsections, to distinguish it from parameters in the reference schemes.

The ETD scheme in for the solution of \eqref{eq...ac} is given as follows.
\begin{align*}
u^{n+1} = \e^{\tau A}u^n + \tau\varphi_1(\tau A)f(u^n),
\end{align*}
where $A = \varepsilon_{AC}^2 D_h$, $\varphi_1(z) = \frac{1-\e^z}{z}$. The time step $\tau$ in this scheme is also referred to as $\tau_{ETD}$. To ensure that the ETD method used in our implementations is state-of-the-art, the standard Krylov method is applied to fast approximate phi-functions as shown in \cite{Goeckler2014}. The training of the neural network was conducted on the Tianhe-2 supercomputer with one GPU core and six CPU cores, while all the remaining parts were performed on a MacBook Air equipped with an Apple M1 chip. 

\subsection{\texorpdfstring{$1$}--dimensional example}\label{sec:1d}
We first provide a $1$-dimensional example where the interfacial width is fixed to $\varepsilon_{AC} = 0.01$. The space interval $\Omega=[-\pi,\pi]\subset\mathbb{R}$ is discretised into $512$ uniform mesh points. The initial data $u_0$ are generated by 
\begin{equation}\label{eq...1dinitdata}
    u_0 = \sum_{i = 1}^{128} a_i \e^{-i/4} \sin(ix)+b_i \e^{-i/4} \cos(ix)
\end{equation}
where the coefficients $a_i$, $b_i \sim \mathcal{N}(0,1)$ are sampled from a normal distribution. We then normalise the generated data by $u_0\leftarrow\frac{u_0}{max(|u_0|)}$ so that $\|u_0\|_{L^\infty}\leq1$. Using this method, we randomly generate $3520$ different initial data and separate them into the training dataset $\{u_0^{(i)}\}_{i=1}^{N=3200}$ and testing dataset $\{u_0^{(i)}\}_{i=1}^{320}$. 

The architecture of the neural network is illustrated in Section \ref{sec:hybrid_methods}. Here we give more details on the hyperparameter setting in this case. The whole neural network is composed of $8$ convolutional layers, each consisting of one Conv1d operator and the activation function $\sigma=Tanh()$. All the $Conv1d$ operators in this case have $21$ kernel size, $1$ stride, and $10$ padding with the \textit{reflect} padding mode. Under this setting, a simplified convolutional neural network with $113409$ learnable parameters is applied as our model to approximate the implicit midpoint time-stepper \eqref{scheme...midpoint} with $\tau_{midpoint}=0.5,1,2$, respectively.
\begin{remark}\label{rmk:choice_of_convolutional_kernel}
    Since the choice of the convolutional kernel size is tricky in practice, we hereby provide a remark to clarify the effect of the kernel size. Note that one way to think of the role of a convolutional kernel is to compare it to a stencil in finite difference schemes. The stencil size typically does not need to change with mesh size, and in a similar vein we found that changing the size of the convolutional kernel is not necessary to maintain similar levels of accuracy over a range of mesh sizes. Of course, the choice of convolutional kernel size may affect the accuracy of neural network prediction, since it will change the number of parameters in the network. 
    However, in practice we found that this effect is not correlated with the meshsize.
\end{remark}

\begin{figure}[htbp]
	\centering
	\subfigure[]{
		\begin{minipage}[b]{0.45\textwidth}
			\includegraphics[width=1\textwidth]{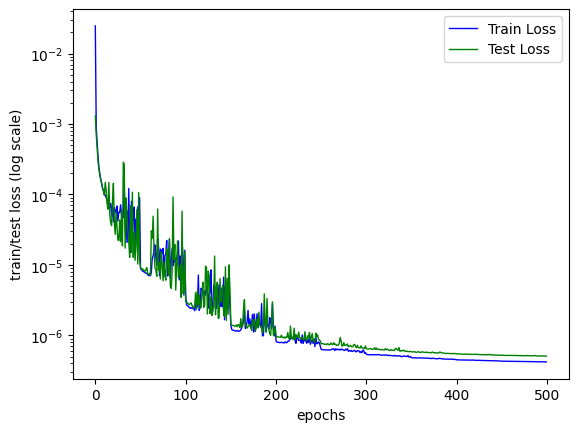}
		\end{minipage}
		\label{fig:1d loss}
	}
    	\subfigure[]{
    		\begin{minipage}[b]{0.45\textwidth}
   		 	\includegraphics[width=1\textwidth]{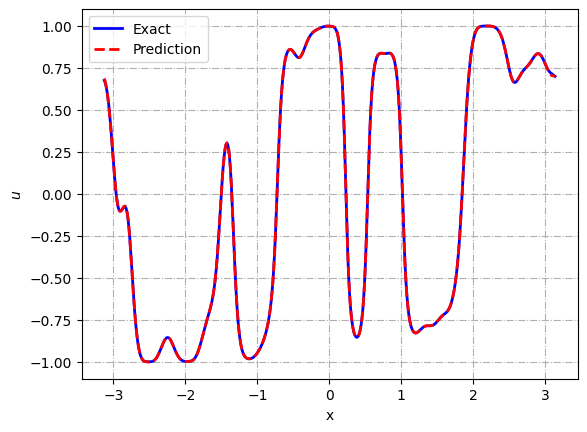}
    		\end{minipage}
		\label{fig:1d test model}
    	}
	\caption{(a) Loss curve; (b) Comparison between neural network prediction and exact midpoint solution}
	\label{fig:1d loss+pred}
\end{figure}

The loss function is defined in \eqref{eq...loss} and is minimized by the optimizer ADAM as mentioned in Section \ref{sec:scheme_informed_learning} with initial learning rate $4\times10^{-4}$ and $L^2$ regularisation weight decay rate $10^{-7}$. The learning rate decays by half after every $50$ epochs and the training process lasts for $500$ epochs. The training process cost around $225$ (sec) for 1D data. We show the training results of $\tau_{midpoint}=1$ case in Figure \ref{fig:1d loss+pred}. Figure \ref{fig:1d loss} shows the loss curve on the training and testing datasets throughout the training process, while \ref{fig:1d test model} compares the neural network prediction with the numerical solution given by the implicit midpoint method with time step $\tau=1$ solved by Newton's method. It is noteworthy that in Figure \ref{fig:1d test model}, the tested initial data is also generated by \eqref{eq...1dinitdata} but lies out of training and testing datasets. The $L^2$ error of the prediction here is lower than $1.67\times10^{-3}$, which indicates the efficiency of the Neural Network approximation to the one step implicit method even with a relatively large time step such as $1$. 

In the following, we use the neural network prediction as the initial guess in Newton's method to solve \eqref{eq...ac} until steady state and compare the performance of our neural hybrid solver with one explicit solver, the ETD scheme, with different time steps and two other Newton solvers -- one with direct initial guess (the numerical solution from the former time step) and the other with ETD solution as the initial guess. In our implementation, all the calculations of the exponential function and $\varphi_1$ function in the ETD scheme are accelerated by the standard Krylov method with generating only $10$ basis vectors in the Arnoldi iteration.
For the simplicity of our formulation, throughout this section we use the following notations to represent the different numerical methods we are considering: 
\begin{itemize}
    \item ``NN hybrid'' refers to the neural network initialised implicit method;
    \item ``ETD hybrid'' refers to the implicit method initialised with the ETD solution;
    \item ``ETD'' refers to the fully explicit ETD scheme.
\end{itemize}

\begin{figure}[htbp]
	\centering
       \subfigure[]{
		\begin{minipage}{0.315\textwidth}
			\includegraphics[width=1\textwidth]{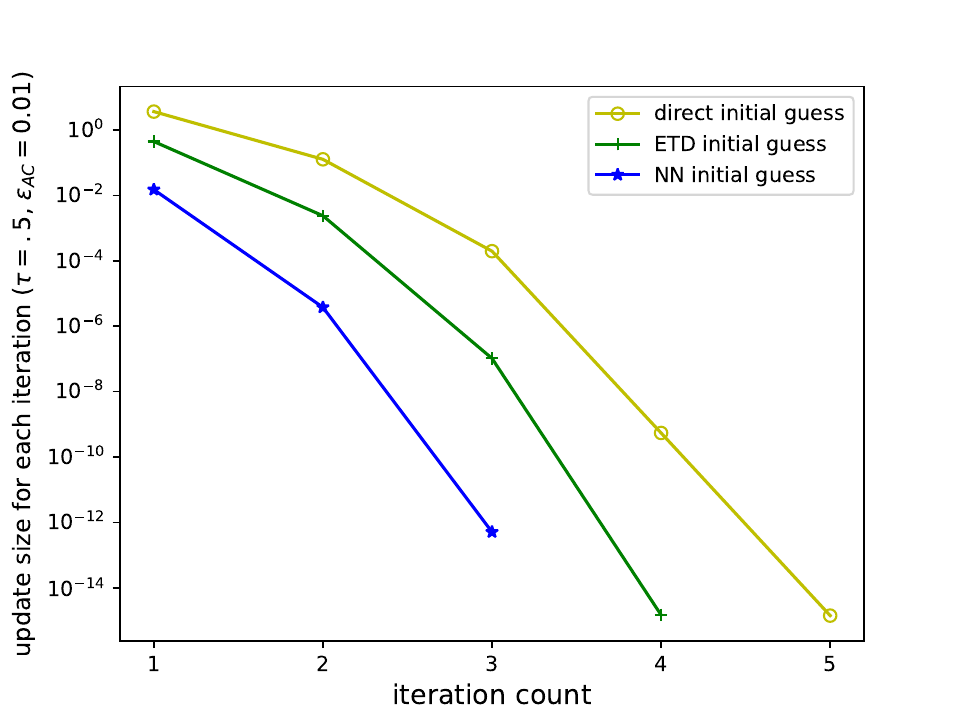}
		\end{minipage}
		\label{fig:1d div_0.5}
	}%
    	\subfigure[]{
    		\begin{minipage}{0.315\textwidth}
   		 	\includegraphics[width=1\textwidth]{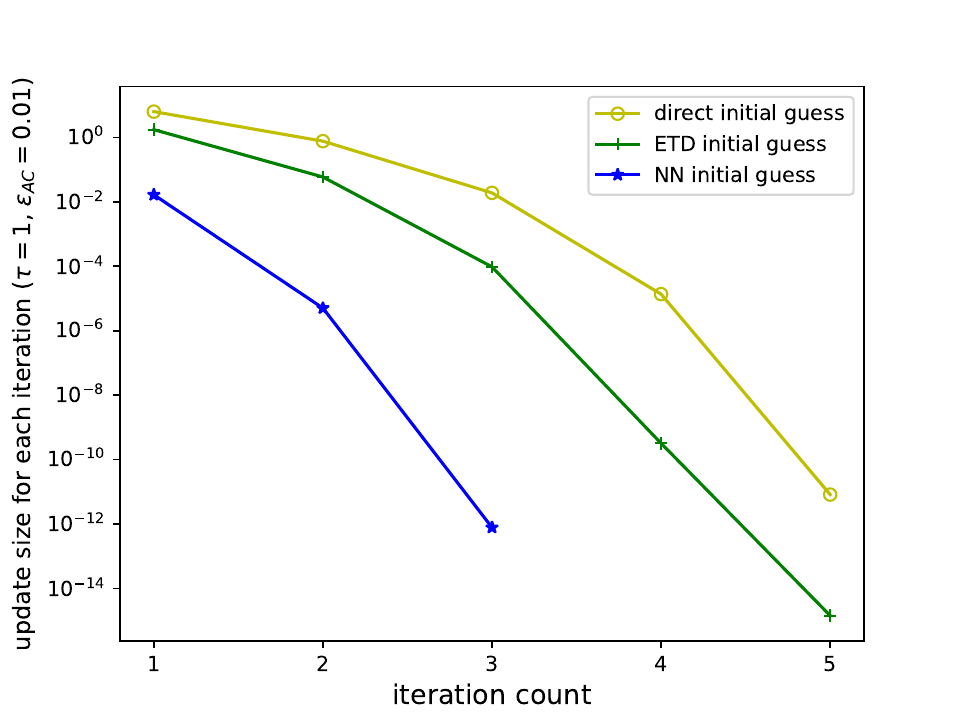}
    		\end{minipage}
		\label{fig:1d div_1}
    	}%
            \subfigure[]{
    		\begin{minipage}{0.315\textwidth}
   		 	\includegraphics[width=1\textwidth]{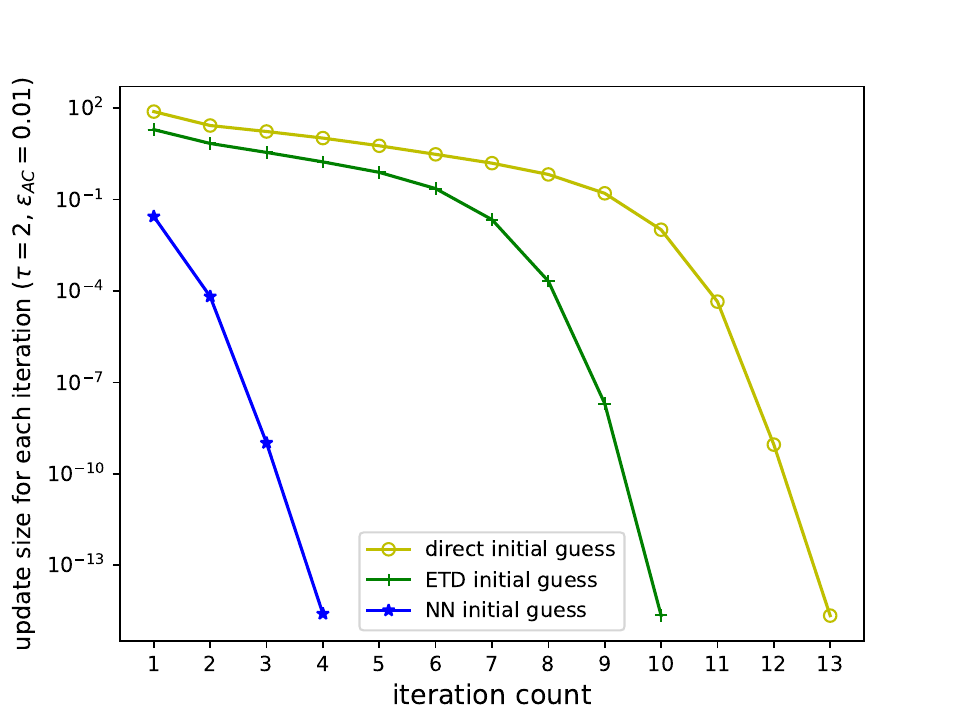}
    		\end{minipage}
		\label{fig:1d div_2}
    	}
	\caption{Update size VS iteration count (1D). The yellow line represents iteration with direct initial guess, the blue line represents iteration with NN prediction initial, and the green line represents iteration with ETD solution initial. (a) $\tau_{midpoint}=0.5$; (b) $\tau_{midpoint}=1$; (b) $\tau_{midpoint}=2$.}
	\label{fig:1d div in iteration}
\end{figure}

To begin with, we check the convergence and the iteration count of Newton iteration using three different initial guesses in solving one step of the implicit midpoint method. Figure \ref{fig:1d div in iteration} exhibits how the update size changes in the iteration. As we can see from the plot, the neural network prediction can provide an initial guess that is much closer to the solution compared to the direct initial guess and the ETD solution with the same time step, which leads to a reduction in iteration count. Furthermore, to eliminate the impact of randomness in initial values, we randomly generate $100$ different initial data and estimate the average iteration count. The result is shown in Table \ref{table: onestep_itercount}.

\begin{table}[htbp]
\scriptsize
\centering
\begin{tabular}{|c|c|c|c|}
\hline
 \diagbox{$\tau_{midpoint}$}{Iteration Count}{Initial Guess}   & Direct initial & ETD ($\tau=\tau_{midpoint}$) & NN prediction \\ \hline
 0.5   & 5.00              & 4.00              & 3.02          \\ \hline
 1   & 5.00              & 5.00              & 3.18          \\ \hline
 2   & 12.10              & 9.57              & 4.25          \\ \hline
\end{tabular}
\caption{Average iteration count in one time step using different initial guess ($1$-dimensional)}
\label{table: onestep_itercount}
\end{table}

\begin{figure}[htbp]
	\centering
	\subfigure[]{
		\begin{minipage}[a]{0.495\textwidth}
			\includegraphics[width=1\textwidth]{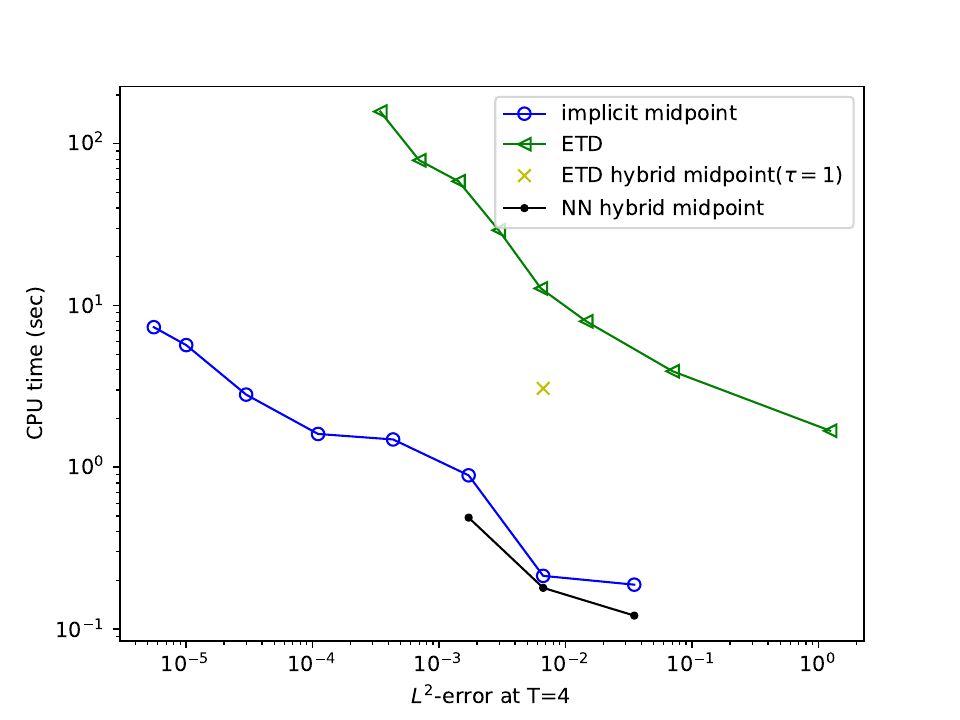}
		\end{minipage}
		\label{fig:1d CPUerr}
	}
 
    	\subfigure[]{
    		\begin{minipage}[b]{0.495\textwidth}
   		 	\includegraphics[width=1\textwidth]{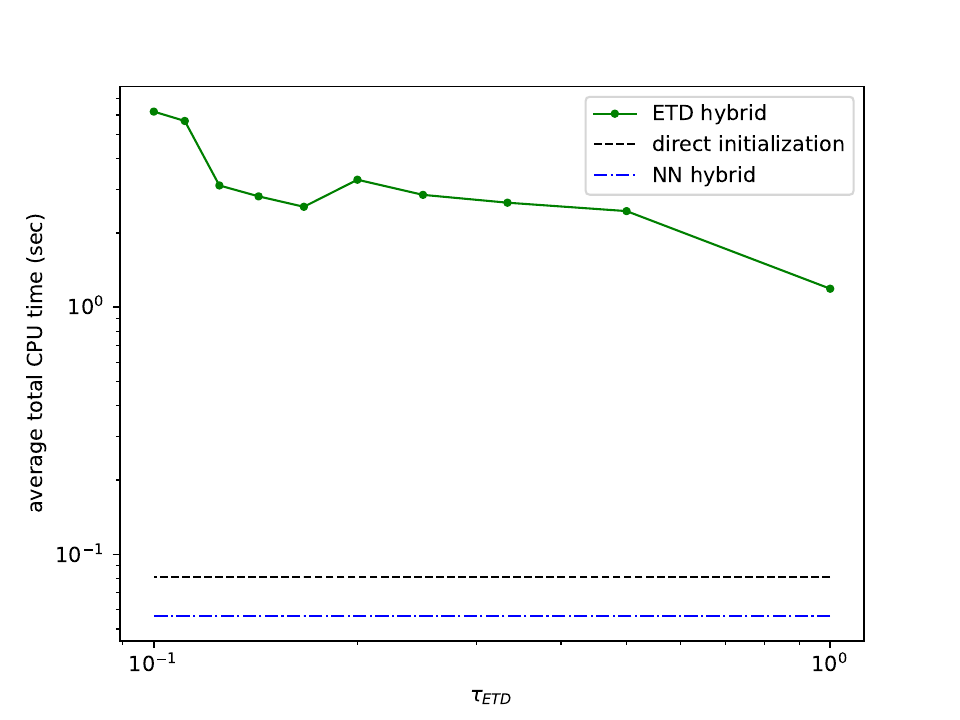}
    		\end{minipage}
		\label{fig:1d singstep totalcost}
    	}%
            \subfigure[]{
    		\begin{minipage}[b]{0.495\textwidth}
   		 	\includegraphics[width=1\textwidth]{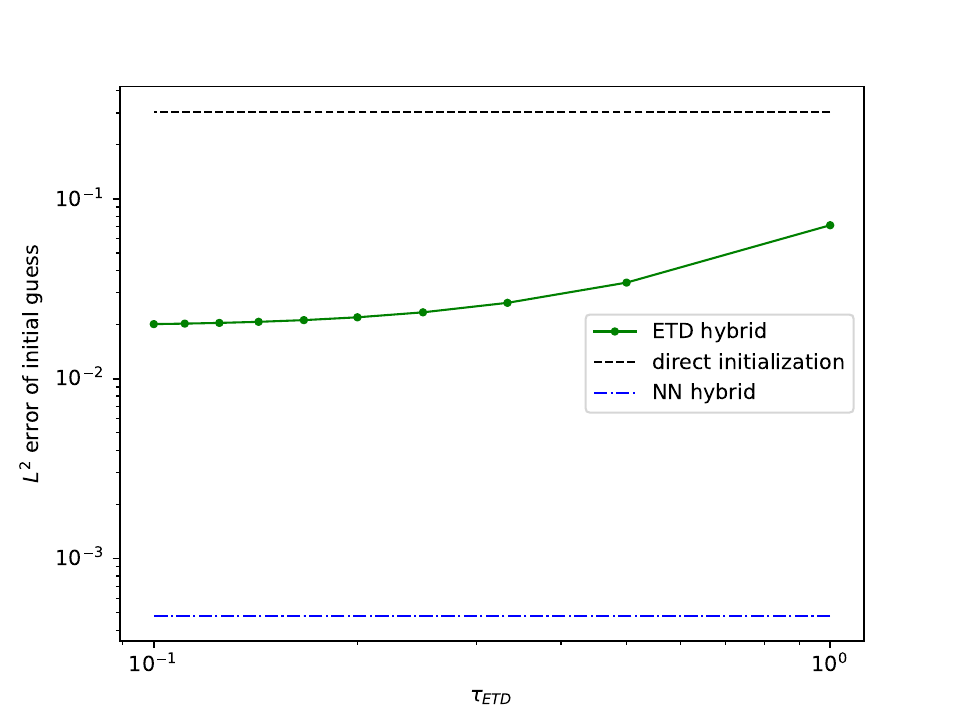}
    		\end{minipage}
		\label{fig:1d singstep initialerror}
    	}
	\caption{(a) CPU time VS the $L^2$ error. The green line represents the performance of ETD scheme, the blue line represents the midpoint method, the black line and the yellow cross are the result of NN hybrid and ETD hybrid solvers, respectively; (b) average CPU time in a single time step when ETD solutions with different time steps are used as the initial guess. The green line shows the average CPU time of ETD scheme, the black and blue lines represent the CPU time of direct initial guess and the NN hybrid method, respectively. Here, we fix $\tau_{midpoint}=1$. (c) $L^2$ error of initial guesses when $\tau_{ETD}$ varies from $0.1$ to $1$. Here, $\tau_{midpoint}=1$.}
	\label{fig:1d CPU+singstepcost}
\end{figure}

Then we consider multiple time steps of the implicit midpoint method and solve the equation until $T=4$. We first plot the CPU time versus the $L^2$ error of solutions at $T=4$ using the implicit midpoint method and the ETD scheme in Figure \ref{fig:1d CPUerr}. To achieve the same accuracy, the implicit midpoint method requires less computational time than the ETD scheme. When $\tau_{midpoint}=1$, the $L^2$ error of the midpoint method solution using direct initial guess is $6.67\times10^{-3}$ and the corresponding CPU time is $0.21$ (sec). The neural hybrid solver reduces the CPU time to $0.18$ (sec) without changing the accuracy, which accelerates the whole computation by $14.29\%$. The same results are observed in the cases when $\tau_{midpoint}=0.5\text{ and }2$ where our neural hybrid solver accelerates the midpoint method by $44.94\%$ and $36.84\%$, respectively. We summarise the result of CPU time of direct initial solver and neural hybrid solver in Table \ref{table: CPUtime} which can be found on page \pageref{table: CPUtime} in order to allow for comparison to the 2-dimensional case. ETD scheme, although it can also produce a better initial guess compared to the direct initial guess, fails to speed up the algorithm because of the high computational cost of its own. To further investigate whether the ETD solution with smaller time steps is capable of decreasing iteration count and the CPU time, we fix $\tau_{midpoint}=1$ and consider a single time step midpoint method equipped with three Newton solvers. We vary $\tau_{ETD}$ in the ETD hybrid solver from $1$ to $0.1$. The green line with dot-markers in Figure \ref{fig:1d singstep totalcost} shows the CPU time (average results of 100 different initial data) of the ETD hybrid solver. The two benchmarks are CPU time using direct initial guess (black line) and the neural network predicted initial guess (blue line). Note that in this single time step, the CPU time for direct initial guess is $8.13\times10^{-2}$ (sec), while the neural hybrid solver costs only $5.66\times10^{-2}$ (sec) which speeds up one step of the implicit method by $30.38\%$. However, the CPU time for the ETD hybrid solver is always higher than the two benchmarks. Figure \ref{fig:1d singstep initialerror} exhibits the $L^2$ error of the three types of initial guesses. Combining the analysis in Figure \ref{fig:1d singstep totalcost} and Figure \ref{fig:1d singstep initialerror}, a straightforward observation here is that if we keep decreasing $\tau_{ETD}$ in the ETD hybrid solver, the CPU time will increase since the calculation of initial guess should be more and more expensive. 

\begin{figure}[htbp]
	\centering
	\subfigure[]{
		\begin{minipage}[b]{0.45\textwidth}
			\includegraphics[width=1\textwidth]{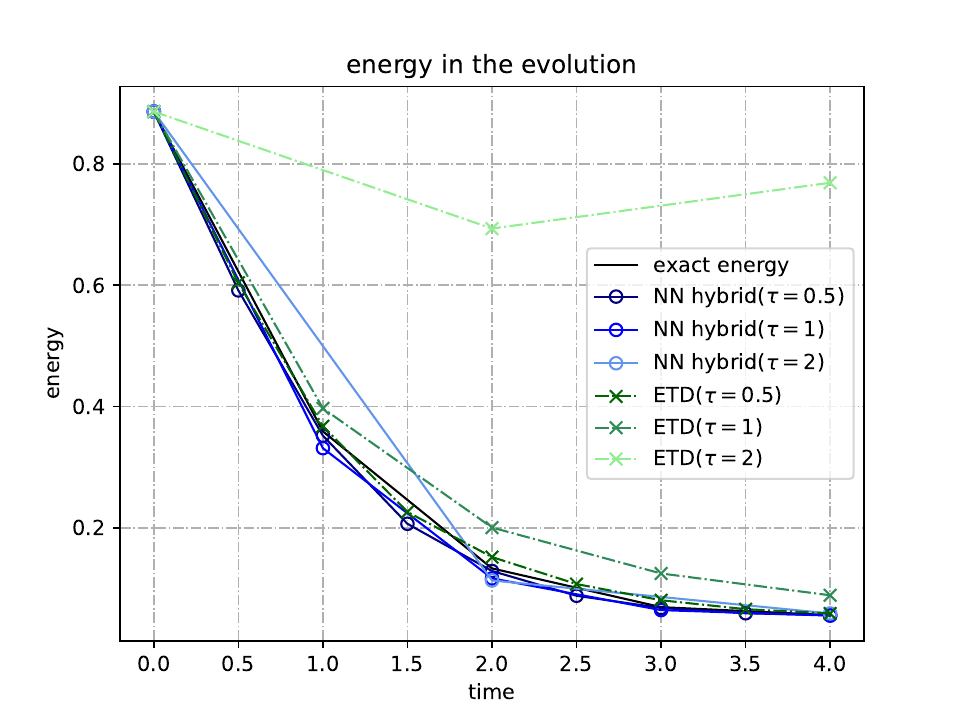}
		\end{minipage}
		\label{fig:1d energy}
	}
    	\subfigure[]{
    		\begin{minipage}[b]{0.45\textwidth}
   		 	\includegraphics[width=1\textwidth]{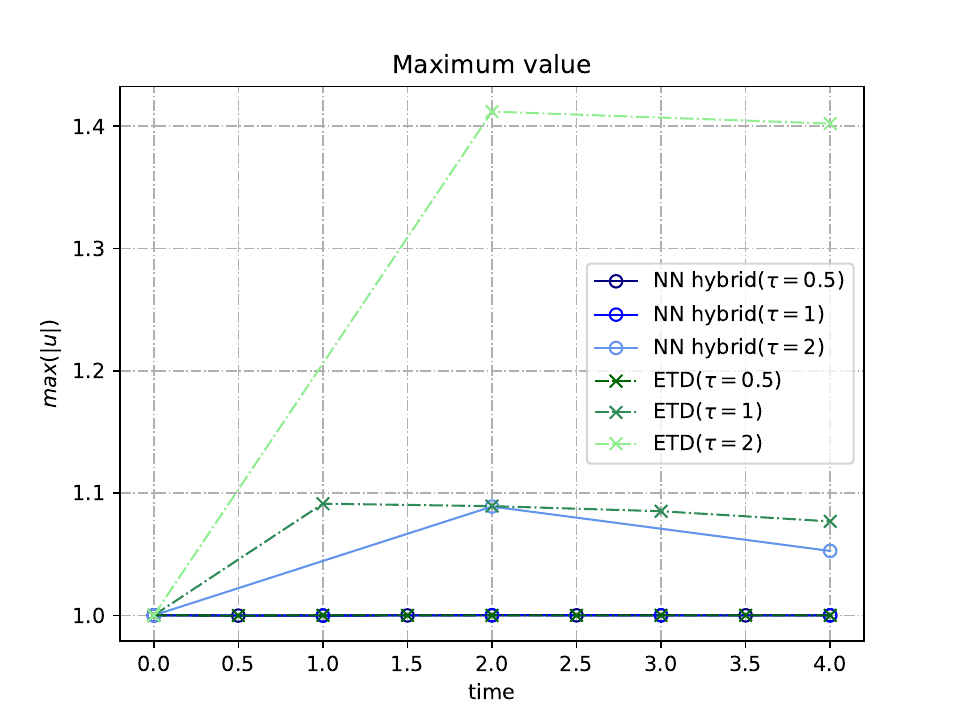}
    		\end{minipage}
		\label{fig:1d maxvalue}
    	}
	\caption{(a) energy of the solution calculated by ETD and neural hybrid solver; (b) maximum absolute value of the solution calculated by ETD and neural hybrid solver.}
	\label{fig:1d structure preserve}
\end{figure}

Finally, we investigate in Figure \ref{fig:1d structure preserve} the structure-preserving properties of the proposed neural hybrid solver. As shown in Figure \ref{fig:1d energy}, the ETD scheme with $\tau_{midpoint}=2$ is unable to preserve the energy dissipation of the solution, whereas the hybrid solver with $\tau_{midpoint}=0.5, 1, 2$ are all able to preserve this property. The maximum absolute value of solution $max(|u|)$ is plotted in Figure \ref{fig:1d maxvalue}. We observe that the ETD scheme with $\tau_{ETD}=1\text{ and }2$ fails to preserve the MBP, while the hybrid solver with $\tau_{midpoint}=0.5\text{ and } 1$ appears to perfectly preserve this property until $T=4$. Even though the hybrid solver with time step $2$ also fails to bound the solution in $[-1,1]$, it still outperforms the ETD scheme with $\tau_{ETD}=1$. 

\subsubsection{Impact of the neural network architecture and mesh resolution on predictive performance}\label{sec:impact}
Following the above experiments, and the comments in Section \ref{sec:nn_architect}, we provide a brief discussion and further evidence supporting our choice of neural network architecture. Here we first present the results of an initial experiment where we replace our simplified convolutional neural network with a fully-connected neural network (FCNN) with $30423552$ parameters. We take the time step of midpoint method to be $1$ and all the hyperparameters in the training are consistent with the experiments using CNN. Table~\ref{table:FCN} shows the average iteration count and the average CPU time required for Newton's method with $100$ random initial data when using four different initial guesses: direct initial, ETD solution ($\tau_{ETD}=\tau_{midpoint}=1$), FCNN prediction and CNN prediction. The results are consistent with the discussion in Section \ref{sec:nn_architect}. In particular, while the fully connected FCNN architecture is able to reduce iteration count it requires significantly more time for inference, thus leading to less reduction in the overall CPU time of the hybrid method. This experiment also lends additional support to our theoretical analysis of the generalisation error in Section~\ref{sec:generalisation_error}: since both architectures were trained with scheme-informed loss \eqref{eq...loss}, we expect their performance on unseen data to be good so long as the Lipschitz constants of the trained neural networks are not too large. This justifies the choice of CNN architecture over a fully connected neural network.
\begin{table}[h!]
\scriptsize
\centering
\begin{tabular}{|c|c|c|c|c|}
\hline
 Initial Guess   & Direct initial & ETD ($\tau=\tau_{midpoint}$) & FCNN prediction & CNN prediction\\ \hline
 Iteration Count   & 5.00              & 5.00              & 4.0 &    \textbf{3.18}     \\ \hline
 CPU time (sec)   & 0.18              & 1.62              & 0.23 &     \textbf{0.10}    \\ \hline
\end{tabular}
\caption{Average iteration count and CPU time in one time step using different initial guesses.}
\label{table:FCN}
\end{table}

Secondly, in order to explore the impact of the mesh resolution on the performance of the neural hybrid method, we repeat the experiments from Table~\ref{table: onestep_itercount} (for $\tau_{midpoint}=1$) with a finer spatial mesh resolution, namely with $1024$ meshpoints. The results are shown in Table~\ref{table:1024mesh} and we included for completeness the case for 512 meshpoints already presented in the previous section. The results show a comparable trend in both mesh resolutions: the NN prediction is able to significantly reduce the iteration count and total cost thus suggesting that the CNN architecture used leads to speed-up which is not very sensitive to the mesh-resolution (cf. interpretation as computational stencils in Remark~\ref{rmk:choice_of_convolutional_kernel}).
\begin{table}[htbp]
\scriptsize
\centering
\begin{tabular}{|c|c|c|c|c|}
 \hline Mesh resolution &Initial Guess   & Direct initial & ETD ($\tau=\tau_{midpoint}$) & NN prediction \\ \hline
\multirow{2}{*}{512}& Iteration Count   & 5.00              & 5.00              & \textbf{3.18 }         \\ \cline{2-5}
 &CPU time (sec)   & 0.08              & 1.19              & \textbf{0.06 }         \\ \hline
\hline
\multirow{2}{*}{1024} &Iteration Count   & 5.00              & 4.00              & \textbf{2.50}          \\ \cline{2-5}
& CPU time (sec)   & 1.22              & 19.19              & \textbf{0.57}          \\ \hline
\end{tabular}
\caption{Average iteration count and CPU time in one time step using different initial guesses ($1$-dimensional).}
\label{table:1024mesh}
\end{table}

\subsection{\texorpdfstring{$2$}--dimensional example}\label{sec:2d}
We now evaluate the proposed method's ability to deal with $2$-dimensional problems. The interfacial width $\varepsilon_{AC}$ is set to $0.02$. The $128\times128$ uniform mesh is applied to discretise the 2D space domain $\Omega=[-\pi,\pi]^2\subset\mathbb{R}^2$. The initial data $u_0$ in this example are generated as follows.
\begin{align*}
u_0(x_1,x_2) = &\sum_{i=1}^{M_1} \sum_{j=1}^{M_2} [a_{ij} \sin(i x_1)\sin(jx_2) + b_{ij} \sin(i x_1)\cos(jx_2)\\
&+ c_{ij} \cos(i x_1)\sin(jx_2)+d_{ij} \cos(i x_1)\cos(jx_2)],\\
\end{align*}
where $a_{ij}, b_{ij},c_{ij},d_{ij}\sim \mathcal{N}(0,1)$ are real coefficients randomly sampled from the normal distribution. Similar to the $1$-dimensional case, all the initial data are normalised such that $u_0\leftarrow\|u_0\|_{L^\infty}\leq1$. We randomly generate $3520$ different initial data and separate them into the training dataset $\{u_0^{(i)}\}_{i=1}^{N=3200}$ and testing dataset $\{u_0^{(i)}\}_{i=1}^{320}$. 

The hyperparameter setting in this example is given as follows. $6$ convolutional layers with one $Conv2d$ operator and the activation function $\sigma=Tanh()$ inside each layer are applied in our neural network model. All the $Conv2d$
operators are equipped with the kernel size $(9,9)$, stride size $(1,1)$, padding size $(4,4)$ and the \textit{reflect} padding mode. With this setting, the entire network has $417473$ learnable parameters.

\begin{figure}[htbp]
	\centering
	\subfigure[]{
		\begin{minipage}[b]{0.45\textwidth}
			\includegraphics[width=1\textwidth]{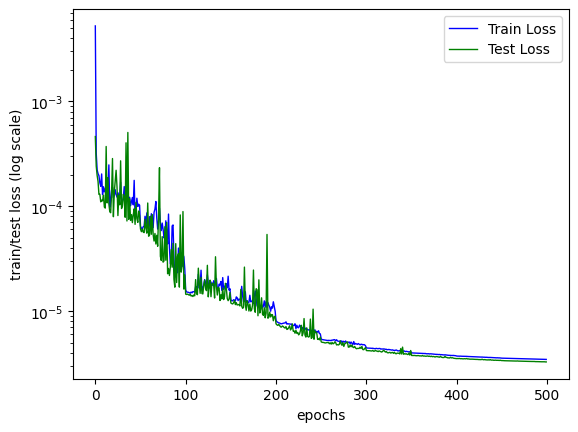}
		\end{minipage}
		\label{fig:2d loss}
	}
 
    	\subfigure[]{
    		\begin{minipage}[b]{0.65\textwidth}
   		 	\includegraphics[width=1\textwidth]{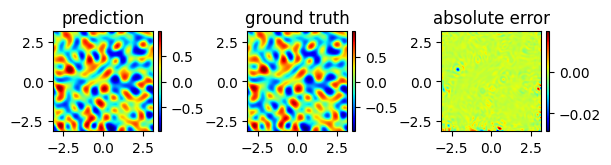}
    		\end{minipage}
		\label{fig:2d test model}
    	}
	\caption{(a) Loss curve; (b) comparison between prediction and exact midpoint solution}
	\label{fig:2d loss+pred}
\end{figure}

The constructed network is then trained by ADAM optimizer to approximate the implicit midpoint time-stepper with $\tau_{midpoint}=0.25, 0,5, 1$. The learning rate starts from $4\times10^{-4}$ and decays by half after every $50$ epochs until $500$ epochs. The $L^2$ regularisation weight decay rate is fixed to $10^{-7}$. The training process cost around $4320$ (sec) for 2D data. We first visualize the training result of the above-described neural network. Figure \ref{fig:2d loss+pred} displays the performance of neural time-stepper with $\tau_{midpoint}=0.5$. Figure \ref{fig:2d loss} shows the training and testing loss throughout $500$ epochs, and Figure \ref{fig:2d test model} compares the network prediction to the exact solution of the implicit midpoint stepper with $\tau_{midpoint}=0.5$. The trained neural time-stepper achieves an $L^2$ absolute error lower than $1.90\times10^{-3}$.

The neural hybrid solver is then investigated in solving this $2$-dimensional problem. We still compare its behaviour with the explicit ETD scheme. With the help of the standard Krylov method, the size of discrete function value vector is reduced from $128^2$ to $16$ in approximating phi-functions, which vastly alleviates the cost of the ETD method.

\begin{figure}[htbp]
	\centering
       \subfigure[]{
		\begin{minipage}{0.3\textwidth}
			\includegraphics[width=1\textwidth]{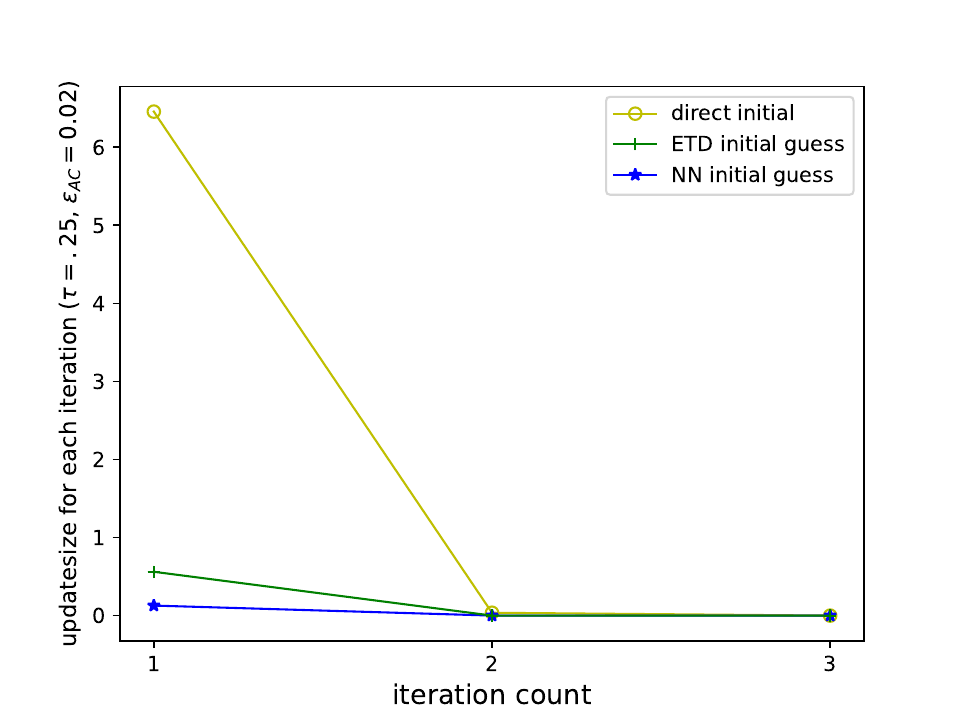}
		\end{minipage}
		\label{fig:2d div_0.25}
	}
    	\subfigure[]{
    		\begin{minipage}{0.3\textwidth}
   		 	\includegraphics[width=1\textwidth]{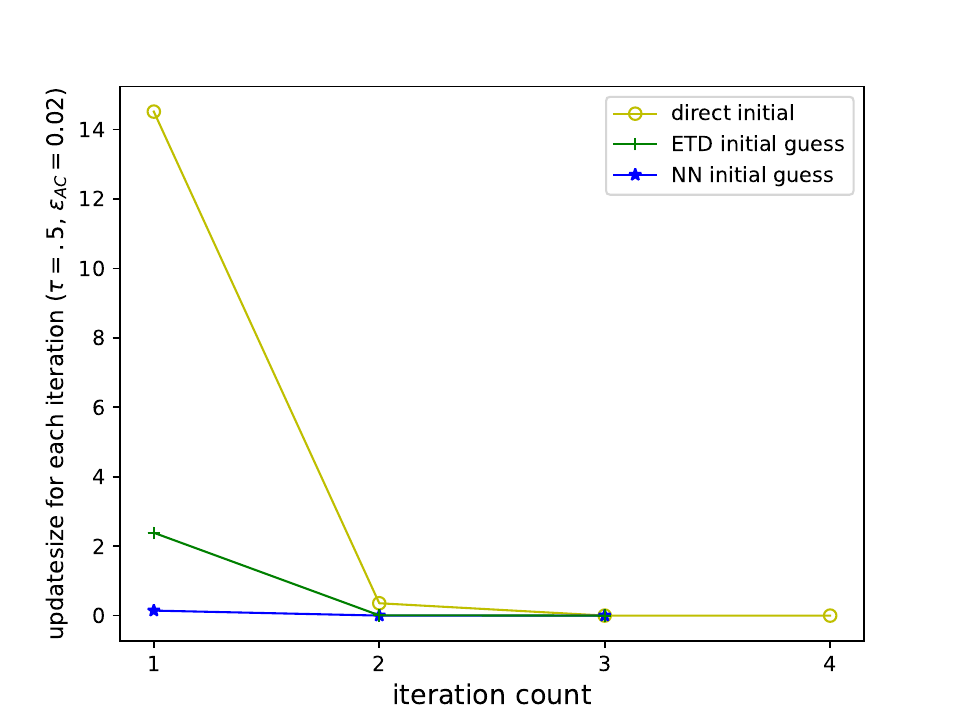}
    		\end{minipage}
		\label{fig:2d div_0.5}
    	}
            \subfigure[]{
    		\begin{minipage}{0.3\textwidth}
   		 	\includegraphics[width=1\textwidth]{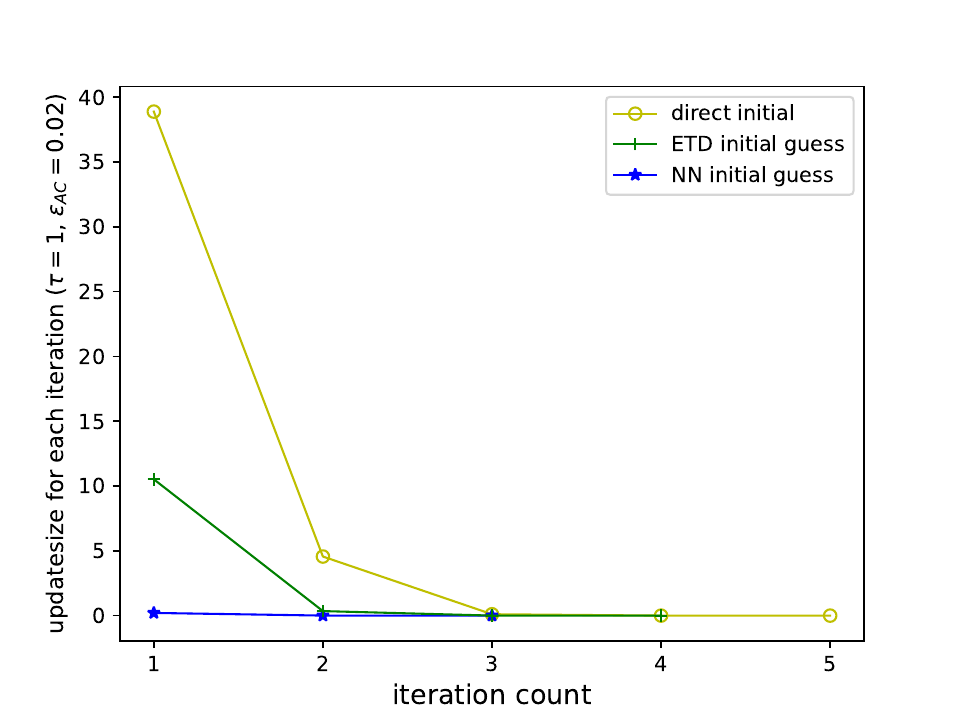}
    		\end{minipage}
		\label{fig:2d div_1}
    	}
	\caption{Update size VS iteration count (2D). The yellow line represents iteration with direct initial guess, the blue line represents iteration with NN prediction initial, and the green line represents iteration with ETD solution initial. (a) $\tau_{midpoint}=0.25$; (b) $\tau_{midpoint}=0.5$; (b) $\tau_{midpoint}=1$.}
	\label{fig:2d div in iteration}
\end{figure}

\begin{figure}[htbp]
	\centering
	\subfigure[]{
		\begin{minipage}[a]{0.45\textwidth}
			\includegraphics[width=1\textwidth]{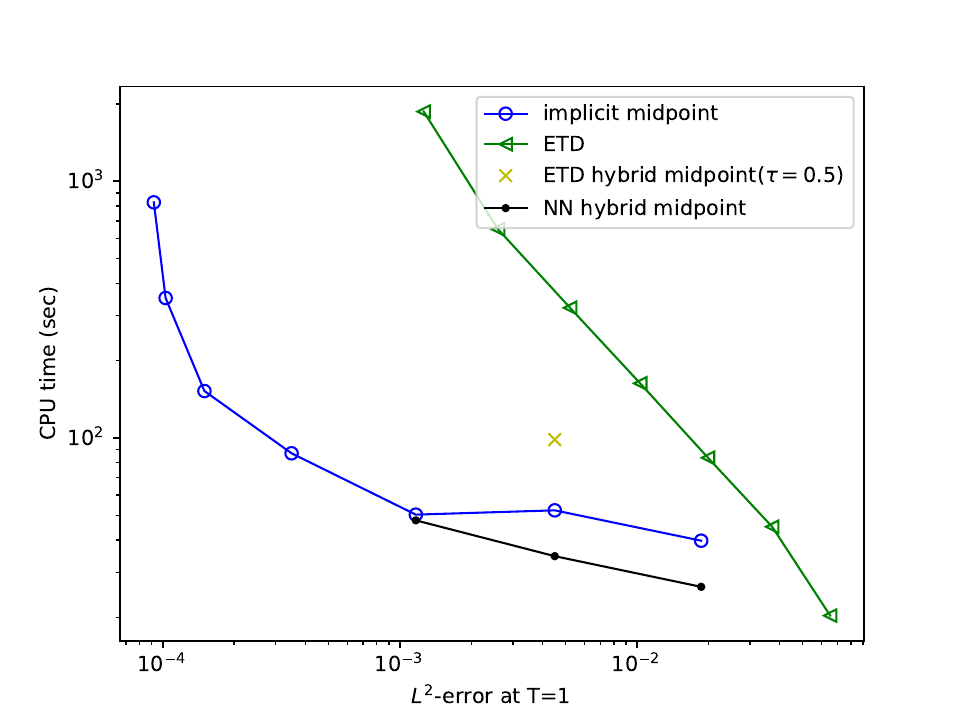}
		\end{minipage}
		\label{fig:2d CPUerr}
	}
 
    	\subfigure[]{
    		\begin{minipage}[b]{0.45\textwidth}
   		 	\includegraphics[width=1\textwidth]{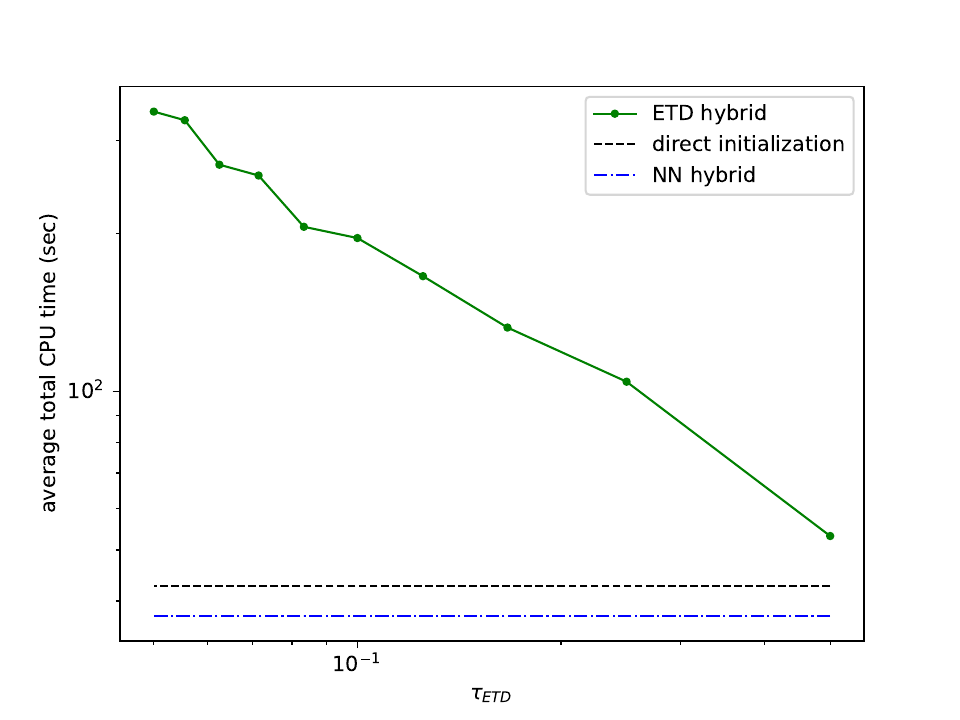}
    		\end{minipage}
		\label{fig:2d singstep totalcost}
    	}
           \subfigure[]{
    		\begin{minipage}[b]{0.45\textwidth}
   		 	\includegraphics[width=1\textwidth]{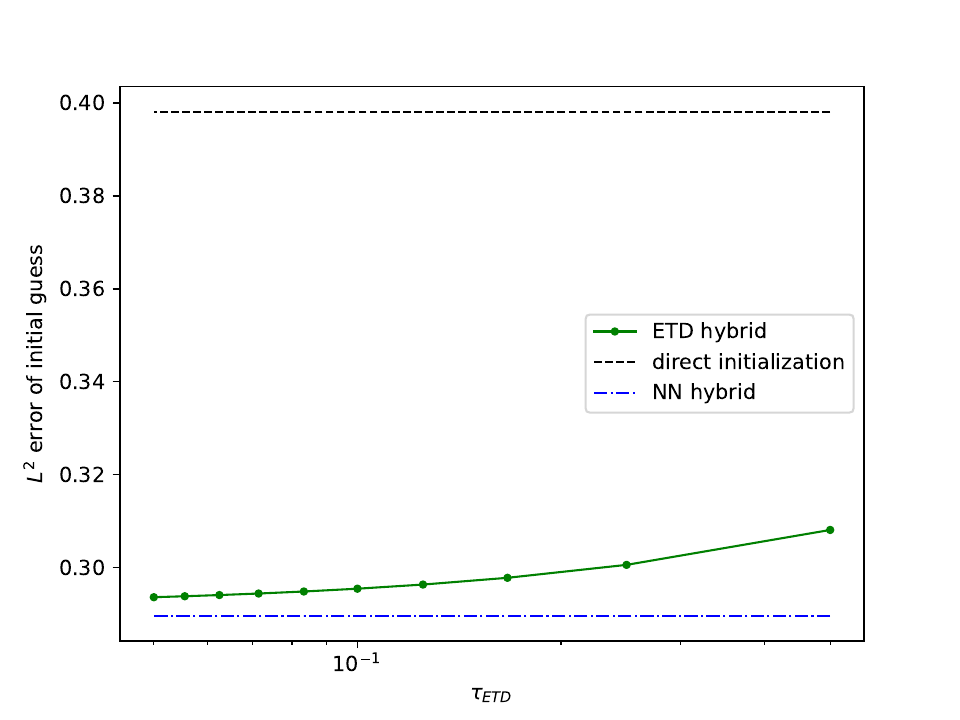}
    		\end{minipage}
		\label{fig:2d singstep initialerror}
    	}
	\caption{(a) CPU time VS the $L^2$ error. The green line represents the performance of the ETD scheme, the blue line represents the midpoint method, the black line and the yellow cross are the result of NN hybrid and ETD hybrid solvers, respectively; (b) average CPU time in a single time step when ETD solutions with different time steps are used as the initial guess. The green line shows the average CPU time of the ETD scheme, the black and blue lines represent the CPU time of direct initial guess and the NN hybrid solver, respectively. Here, we fix $\tau_{midpoint}=0.5$; (c) $L^2$ error of initial guess when $\tau_{ETD}$ varies from $0.05$ to $0.5$. Here, $\tau_{midpoint}=0.5$.}
	\label{fig:2d CPU+singstepcost}
\end{figure}

We start again with solving a single step in the implicit midpoint method by Newton's method with three different initial guesses. Figure \ref{fig:2d div in iteration} exhibits the update size in the iteration. Notice that the dissipation of the update size indicates the convergence of the Newton's iteration. As shown in the plot, both the neural network prediction and the ETD solution with $\tau_{ETD}=\tau_{midpoint}$ offer better initial guesses and are capable of reducing iteration count when the time step is relatively large. We also display the result of the average iteration count over $10$ different initial data in Table \ref{table: onestep_itercount2d}. Notice that, when $\tau_{midpoint}=0.25$ the CPU time for direct initial guess, NN initial guess and ETD initial guess are $33.74, 27.64, 48.62$ (sec), respectively. Thus, although the neural hybrid solver fails to reduce the iteration count in this case, it can still accelerate the algorithm by $18.08\%$. Since we use GMRES for the solution of the linear systems in Newton's iteration we suspect that this means the Krylov subspaces generated by the Neural Network initial guess provide a better approximation to the true solution value of the midpoint rule than the Krylov subspaces generated by either alternative initial guess. The multiple-time-steps result in Figure \ref{fig:2d CPUerr} further supports this conclusion.

\begin{table}[htbp]
\scriptsize
\centering
\begin{tabular}{|c|c|c|c|}
\hline
 \diagbox{$\tau_{midpoint}$}{Iteration Count}{Initial Guess}   & Direct initial & ETD ($\tau=\tau_{midpoint}$) & NN prediction \\ \hline
 0.25   & 3.00              & 3.00              & 3.00          \\ \hline
 0.5   & 4.00              & 3.00              & 3.00          \\ \hline
 1   & 5.00              & 4.00              & 3.00          \\ \hline
\end{tabular}
\caption{Average iteration count in one time step using different initial guess ($2$-dimensional)}
\label{table: onestep_itercount2d}
\end{table}

We then examine the performance in multiple time steps of the three implicit time steppers and one explicit time stepper. Figure \ref{fig:2d CPUerr} shows the CPU time required to achieve a given $L^2$ error at $T=1$ for both the implicit midpoint method and the ETD scheme. Generally speaking, the implicit midpoint method is more computationally efficient than the ETD scheme if the same accuracy is required. The black line with dot-markers indicates that our neural hybrid solver accelerates the algorithm by $4.99\%$, $33.70\%$, and $33.90\%$ when we consider time steps $\tau_{midpoint}=0.25,0.5\text{ and }1$. The results of CPU time and acceleration rate are also exhibited in Table \ref{table: CPUtime}. We highlight that the CPU time evaluated in Table \ref{table: CPUtime} is the overall cost of Newton's method which includes the time to obtain the initial guesses. More detailed discussion of the cost of individual components of the neural hybrid solver can be found in Appendix \ref{appendix:wall-clock}. As displayed by the yellow cross in Figure \ref{fig:2d CPUerr}, the ETD hybrid solver increases the CPU time of the implicit midpoint method. To further explore whether the reduction in iteration count can compensate the cost of computing a better initial guess by choosing a smaller $\tau_{ETD}$, we fix $\tau_{midpoint}$ to be $0.5$ and vary $\tau_{ETD}$ from $0.05$ to $0.5$ in the ETD hybrid solver for a single time step. The outcomes in Figure \ref{fig:2d singstep totalcost} and \ref{fig:2d singstep initialerror} show that although decreasing $\tau_{ETD}$ leads to an improvement in initialisation, this advantage is offset by the increase in CPU time required in the ETD hybrid solver.

\begin{table}[htbp]
\centering
\begin{tabular}{|c|c|c|c|c|}
\hline
Test case           & $\tau_{midpoint}$ & Direct initial guess & Neural hybrid & Acc. Rate \\ \hline
\multirow{3}{*}{1D} & 0.5               & 0.89                 & 0.49          & 44.94\%   \\ \cline{2-5} 
                    & 1                 & 0.21                 & 0.18          & 14.29\%   \\ \cline{2-5} 
                    & 2                 & 0.19                 & 0.12          & 36.84\%   \\ \hline
\multirow{3}{*}{2D} & 0.25              & 50.26                & 47.75         & 4.99\%    \\ \cline{2-5} 
                    & 0.5               & 52.23                & 34.63         & 33.70\%   \\ \cline{2-5} 
                    & 1                 & 39.79                & 26.30         & 33.90\%   \\ \hline
\end{tabular}
\caption{CPU time (sec) required for solving equation \eqref{eq...ac} until T. The acceleration rate of neural hybrid solver is shown in the last column. In 1D case, $T=4$; in 2D case, $T=1$.}
\label{table: CPUtime}
\end{table}

\begin{figure}[htbp]
	\centering
	\subfigure[]{
		\begin{minipage}[b]{0.45\textwidth}
			\includegraphics[width=1\textwidth]{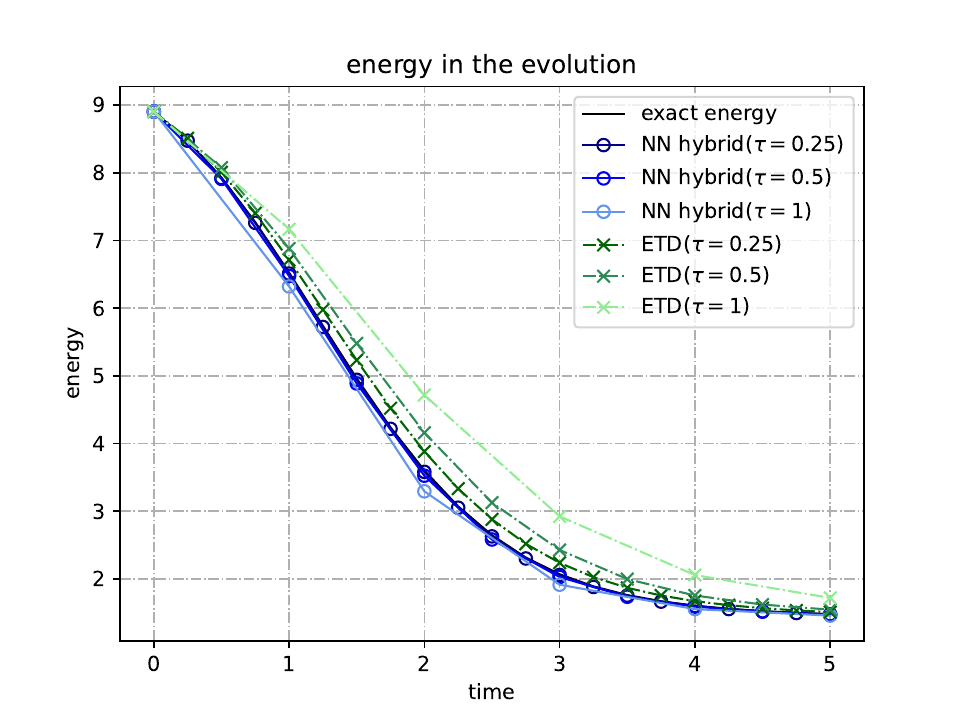}
		\end{minipage}
		\label{fig:2d energy}
	}
    	\subfigure[]{
    		\begin{minipage}[b]{0.45\textwidth}
   		 	\includegraphics[width=1\textwidth]{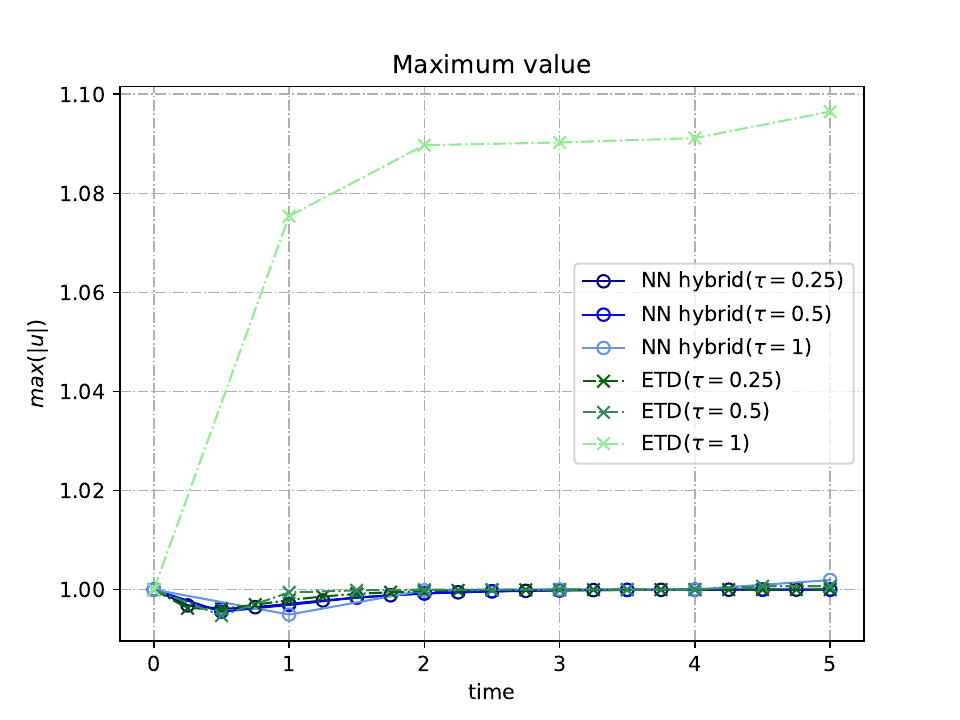}
    		\end{minipage}
		\label{fig:2d maxvalue}
    	}
	\caption{(a) energy of the solution calculated by the ETD and neural hybrid solver; (b) maximum absolute value of the solution calculated by the ETD and neural hybrid solver.}
	\label{fig:2d structure preserve}
\end{figure}

The structure-preserving properties of the neural hybrid solver and ETD scheme are visualized in Figure \ref{fig:2d structure preserve}. Figure \ref{fig:2d energy} shows even though all the solvers can preserve energy dissipation, the neural hybrid solver has a higher energy accuracy than the ETD scheme. We witness in Figure \ref{fig:2d maxvalue} that the ETD scheme with $\tau_{ETD}=1$ fails to preserve the MBP, while all other solvers appear to preserve this property well until $T=5$.

\section{Conclusions and Discussions}\label{sec:conclusion}
In this work, we introduce a neural hybrid solver that accelerates Newton's method for implicit time stepping schemes using improved initialisation provided by a neural network. The neural time-stepper has the structure of a simplified CNN and is trained using a novel implicit-scheme-informed learning strategy. Consequently, it benefits from the advantages of being easy to implement, having fewer parameters leading to fast forward pass and, hence initialisation, and not requiring any labeled data for training. A quantifiable estimate of the reduction in the Newton iteration count for an improved initialisation is provided in Theorem \ref{thm...asym_iteration_count}, which indicates the advantage provided by a good initial guess in Newton's method. To illustrate the robustness of the proposed implicit-scheme-informed learning strategy, we prove in Theorem \ref{thm...genlerror} that the generalisation error can be made arbitrarily close to the training error provided a sufficient amount of training data is used and in Proposition \ref{prop...num_traindata} we provide an upper bound on this number of training data required to achieve a desired accuracy in the generalisation error.

The efficiency of our proposed solver is demonstrated using several applications on solving Allen--Cahn equation equipped with homogeneous Neumann boundary condition using the implicit midpoint method. Numerical implements, both in $1$- and $2$-dimensional settings, show that the neural hybrid solver outperforms other solvers in terms of accuracy, computational cost and structure-preservation. We remark that all the experiments displayed are completed on a pure CPU device after the neural network is trained. We found in our further exploration that given an input data, especially high dimensional data, the time for neural network prediction can be accelerated a lot if it is conducted on GPU devices. This suggests that there is significant potential to extend the neural hybrid solver to a wide range of high dimensional problems in future work. Future research will also investigate how the stiffness of equations affects the training and efficiency of the neural hybrid solver.

Another further direction for future work is the deeper understanding of the number of training data points required in the present set-up. Indeed, while Proposition \ref{prop...num_traindata} provides a reasonable upper bound for this, our numerical experiments indicate that good performance of our methodology is already achieved with a training dataset of more moderate size.

\section*{Acknowledgements}
The authors would like to thank Qinyan Zhou (Sorbonne Universit\'e) for performing initial numerical studies that helped inform this research. Furthermore, the authors express their gratitude to Thibault Faney (IFP Energies Nouvelles), Antoine Lechevallier (IFP Energies Nouvelles \& Sorbonne Universit\'e), and Frédéric Nataf (Sorbonne Universit\'e) for several helpful discussions. T.J. gratefully acknowledges funding from postgraduate studentship and Overseas Research Award of HKUST. K.S. gratefully acknowledges funding from the European Research Council (ERC) under the European Union’s Horizon 2020 research and innovation programme (grant agreement No. 850941). The work of Y.X. was supported by the Project of Hetao Shenzhen-HKUST Innovation Cooperation Zone HZQB-KCZYB-2020083. T.J. and Y.X. would like to thank HKUST Fok Ying Tung Research Institute and National Supercomputing Center in Guangzhou Nansha Sub-center for providing high performance computational resources.
\appendix

\renewcommand{\thesection}{\Alph{section}} 
\makeatletter
\renewcommand\@seccntformat[1]{\appendixname\ \csname the#1\endcsname.\hspace{0.5em}}
\makeatother

\section{Proof of two lemmas}\label{appendix:two_lemma}
We first recall the statement of Lemma \ref{lem:loss_bounds_approximation_error}:
\begin{lemma}
    Suppose $\alpha>\gamma\geq 0$ and $\alpha>d/2$ and that Assumption~\ref{asump...mainthm} is satisfied, and let $G:H^{\gamma}\rightarrow H^{\gamma}$ be a Lipschitz continuous map. Then there is a $\tau_0>0$ and a constant $C>0$ such that for all $u\in B_{H^{\alpha}}(0,1)\cap  B_{L^\infty}(0,1)$ and all $0<\tau\leq \tau_0$ the following estimate holds:
    \begin{align*}
         \|\Phi_{\tau}(u)-G(u)\|_{H^{\gamma}}\leq C\|\mathcal{R}(G(u))\|_{H^\gamma},
\end{align*}
where $\mathcal{R}(G(u)) = \mathcal{R}(G(u);u)=\Psi_\tau(u,G(u))-G(u)$.
\end{lemma}
\begin{proof}
For any $u\in B_{H^\alpha}(0,1)\cap B_{L^\infty}(0,1)$, we have the following identity:
\begin{align}\label{eq...first_identity_implicit_scheme}
    G(u)-\Phi_{\tau}(u)=\Psi_\tau(u,G(u))-\Psi_\tau(u,\Phi_{\tau}(u))-\mathcal{R}(\Phi_{\tau}^*(u))
\end{align}
The implicit midpoint scheme for the Allen--Cahn equation \eqref{eq...ac} gives the following expression for $\Psi_\tau$:
\begin{align*}
    \Psi_\tau(u,v)=u+\tau \varepsilon_{AC}^2\Delta\left(\frac{u+v}{2}\right)-\tau\left(\frac{u+v}{2}\right)\left(\frac{1}{4}(u^2+2uv+v^2)-1\right)
\end{align*}
Using this and \eqref{eq...first_identity_implicit_scheme} we find
\begin{align*}
    G(u)-\Phi_{\tau}(u)&=\frac{\tau}{2} \varepsilon_{AC}^2\Delta\left(G(u)-\Phi_{\tau}(u)\right)+\frac{\tau}{2}\left(G(u)-\Phi_{\tau}(u)\right)\\
    &\quad-\tau\left[\left(\frac{u+G}(u){2}\right)^3-\left(\frac{u+\Phi_\tau(u)}{2}\right)^3\right] \\
    &\quad-\mathcal{R}(\Phi_{\tau}^*(u)).
\end{align*}
Thus
\begin{align}\begin{split}\label{eq...stability_nn_vs_midpt}
    ((1-\frac{\tau}{2})I-\frac{\tau}{2}\varepsilon_{AC}^2\Delta)\left(G(u)-\Phi_{\tau}(u)\right)=
    &-\tau\left[\left(\frac{u+G}(u){2}\right)^3-\left(\frac{u+\Phi_\tau(u)}{2}\right)^3\right]\\
    &-\mathcal{R}(G(u)).
    \end{split}
\end{align}
For $\tau<2$ we have, due to the spectrum of $\Delta$,
\begin{align*}
    C_{1,\tau}:=\|((1-\frac{\tau}{2})I-\frac{\tau}{2}\varepsilon_{AC}^2\Delta)^{-1}\|_{H^\gamma\rightarrow H^\gamma}<\infty.
\end{align*}
Hence we can perform the following estimate:
\begin{align*}
    \|\Phi_{\tau}(u)-G(u)\|_{H^{\gamma}}
    &\leq C_{1,\tau}\frac{\tau}{8}\left\|\left(u+G(u)\right)^3-\left(u+\Phi_\tau(u)\right)^3\right\|_{H^\gamma}\\
    &\quad+C_{1,\tau}\|\mathcal{R}(G(u))\|_{H^\gamma}\\
    &= C_{1,\tau}\frac{\tau}{8}\|[(u+G(u))^2+(u+G(u))(u+\Phi_{\tau}(u))\\
    &\quad+(u+\Phi_\tau(u))^2]\left(G(u)-\Phi_{\tau}(u)\right)\|_{H^\gamma}\\
    &\quad+C_{1,\tau}\|\mathcal{R}(G(u))\|_{H^\gamma}\\
    &\leq \tau C_{2}\|G(u)-\Phi_{\tau}(u)\|_{H^\gamma}+C_{1,\tau}\|\mathcal{R}(G(u))\|_{H^\gamma}
\end{align*}
where Lemma \ref{lemma...sobmulti} is used to obtain the last inequality and 
\begin{align*}
    C_2:=\frac{1}{8}C_{1,\tau}C_\gamma\|(u+G(u))^2+(u+G(u))(u+\Phi_{\tau}(u))+(u+\Phi_\tau(u))^2\|_{H^\gamma}
\end{align*}
Since $G$ is Lipschitz, we can bound $\|G(u)\|_{H^\gamma}$ in terms of (a smooth function of) $\|u\|_{H^\alpha}$. Similarly to the argument in \eqref{eq...stability_nn_vs_midpt} we also have
\begin{align*}
    \Phi_{\tau}(u)=\frac{(1+\frac{\tau}{2})I+\frac{\tau\varepsilon_{AC}^2}\Delta}{2}{(1-\frac{\tau}{2})I-\frac{\tau\varepsilon_{AC}^2}\Delta}{2}u-\frac{\tau}{(1-\frac{\tau}{2})I-\frac{\tau\varepsilon_{AC}^2}\Delta}{2}\left(\frac{u+\Phi_{\tau}(u)}{2}\right)^3.
\end{align*}
Thus, by Assumption \ref{asump...mainthm} we know that for $\tau$ sufficiently small $\|\Phi_{\tau}(u)\|_{H^\gamma}$ is bounded in terms of a Lipschitz function of $\|u\|_{H^\gamma}$.
Moreover, the definition of Sobolev norm implies that $\|u\|_{H^\gamma}\leq \|u\|_{H^\alpha}\leq 1$. 
Thus, $C_2$ is uniformly bounded for any $u\in B_{H^\alpha}(0,1)$ since $\|G(u)\|_{H^\gamma}$, $\|\Phi_{\tau}(u)\|_{H^\gamma}$ and $\|u\|_{H^\gamma}$ are all bounded.
\noindent
We continue to estimate $\|\Phi_{\tau}(u)-G(u)\|_{H^{\gamma}}$. For $\tau<\frac{1}{C_2}$ we have
\begin{align*}
    \|\Phi_{\tau}(u)-G(u)\|_{H^{\gamma}}\leq \frac{C_{1,\tau}}{1-\tau C_2}\|\mathcal{R}(G(u))\|_{H^\gamma}.
\end{align*}
The important thing to recall here is that
\begin{align*}
    \|\mathcal{R}(G(u))\|_{H^\gamma}=\|\Psi_\tau(u,G(u))-G(u)\|_{H^\gamma}.
\end{align*}
In conclusion, let $\tau_0=\min\{2,\frac{1}{C_2}\}$ and $C=\frac{C_{1,\tau}}{1-\tau C_2}$. Then we have the following estimate valid for any $0<\tau<\tau_0$:
\begin{align*}
         \|\Phi_{\tau}(u)-G(u)\|_{H^{\gamma}}\leq C\|\mathcal{R}(G(u))\|_{H^\gamma},
\end{align*}
This completes the proof of this statement.
\end{proof}

We then recall the statement of Lemma \ref{lem:estimate_onresidual}:
\begin{lemma}
Suppose $\alpha> \beta+2>d/2$ and that Assumption \ref{asump...mainthm} is satisfied. Let $G:H^{\beta}\rightarrow H^\beta$ be a Lipschitz continuous map, then there are constants $C,\tau_0>0$ such that for all $0\leq\tau\leq \tau_0$ and any $v,w\in B_{H^\alpha}(0,1)$ we have
\begin{align}
    \|\mathcal{R}(G(v))-\mathcal{R}(G(w))\|_{H^\beta}\leq C\|v-w\|_{H^{\beta+2}}.
\end{align}
\end{lemma}
\begin{proof}
    We have
\begin{align*}
 \|\mathcal{R}(G(v))-\mathcal{R}(G(w))\|_{H^\beta}&=\|\Psi_\tau(v,G(v))-G(v)\\
 &\quad-\Psi_\tau(w,G(w))+G(w)\|_{H^\beta}\\ 
 &\leq \|\Psi_\tau(v,G(v))-\Psi_\tau(w,G(w))\|_{H^\beta}\\
 &\quad+\|G(v)-G(w)\|_{H^\beta}\\
 &\stackrel{\text{def}}{=}\text{\RomanNumeralCaps{1}}+\text{\RomanNumeralCaps{2}}
\end{align*}
Since $F$ is Lipschitz continuous, there is a constant $C_{lip}>0$ such that
\begin{align*}
    \text{\RomanNumeralCaps{2}} = \|G(v)-G(w)\|_{H^\beta}
                               \leq C_{lip}\|v-w\|_{H^\beta}.
\end{align*} 
It remains to estimate \RomanNumeralCaps{1}. Notice that 
\begin{align*}
    \text{\RomanNumeralCaps{1}} &= \|\Psi_\tau(v,G(v))-\Psi_\tau(w,G(w))\|_{H^\beta}\\\
    &=\|((1+\frac{\tau}{2})I+\frac{\tau}{2}\varepsilon_{AC}^2\Delta)(v-w)+(\frac{\tau}{2}I+\frac{\tau}{2}\varepsilon_{AC}^2\Delta)(G(v)-G(w))\\
    &\quad-\frac{\tau}{8}[(v+G(v))^2+(v+G(v))(w+G(w))+(w+G(w))^2]\\
    &\quad\quad\cdot(v-w+G(v)-G(w))\|_{H^\beta}\\
    &\leq \|((1+\frac{\tau}{2})I+\frac{\tau}{2}\varepsilon_{AC}^2\Delta)(v-w)\|_{H^\beta}        +\|(\frac{\tau}{2}I+\frac{\tau}{2}\varepsilon_{AC}^2\Delta)(G(v)-G(w))\|_{H^\beta}\\
    &\quad+\frac{\tau}{8}C_\beta\|(v+G(v))^2+(v+G(v))(w+G(w))+(w+G(w))^2\|_{H^\beta}\\
    &\quad\quad\cdot\|v-w+G(v)-G(w)\|_{H^\beta}\\
\end{align*} 
Lemma \ref{lemma...sobmulti} is used to obtain the last inequality. 
We now note that
\begin{align*}
    \|\Delta u\|_{H^\beta}\leq \|u\|_{H^{\beta+2}}\leq \|u\|_{H^\alpha}.
\end{align*}
Also, from the fact that $\|v\|_{H^\beta}$,$\|w\|_{H^\beta}$,$\|G(v)\|_{H^\beta}$ and $\|G(v)\|_{H^\beta}$ are all bounded, we can introduce another constant $C_3$ which is uniformly bounded for $w\in B_{H^\alpha}(0,1)$:
\begin{align*}
    C_3 = \|(v+G(v))^2+(v+G(v))(w+G(w))+(w+G(w))^2\|_{H^\beta}.
\end{align*}
Substituting this into the estimate of \RomanNumeralCaps{1} yields
\begin{align*}
    \text{\RomanNumeralCaps{1}}&\leq (1+\frac{\tau}{2}+\frac{\tau}{2}\varepsilon_{AC}^2)\|v-w\|_{H^{\beta+2}}
                                +(\frac{\tau}{2}+\frac{\tau}{2}\varepsilon_{AC}^2)C_{lip}\|v-w\|_{H^{\beta+2}}\\
                               &\quad+\frac{\tau}{8}C_\beta C_3(1+C_{lip})\|v-w\|_{H^\beta}\\
                               &\leq C_4 \|v-w\|_{H^{\beta+2}},\\
\end{align*}
where $C_4=1+\frac{\tau}{2}+\frac{\tau}{2}\varepsilon_{AC}^2+(\frac{\tau}{2}+\frac{\tau}{2}\varepsilon_{AC}^2)C_{lip}+\frac{\tau}{8}C_\beta C_3(1+C_{lip})$.
Therefore, the estimates of \RomanNumeralCaps{1} and \RomanNumeralCaps{2} lead to
\begin{align*}
    \|\mathcal{R}(G(v))-\mathcal{R}(G(w))\|_{H^\beta}\leq C\|v-w\|_{H^{\beta+2}},
\end{align*}
where $C=C_{lip}+C_4$.
\end{proof}

\section{Proof of Proposition \ref{prop...num_traindata}}\label{appendix:proof_prop}
Following recall of the statement, the proof of Proposition \ref{prop...num_traindata} is provided.
\begin{proposition}(Proposition \ref{prop...num_traindata})
Fix $\alpha>\gamma>d/2.$ Given $\epsilon>0$ there is an $\epsilon$-Net of size at most
\begin{align}\label{eqn:finite_dim_eps_net}
    N_{\epsilon}^{(M,d)}=\left\lceil\left(\frac{2}{\epsilon}\right)^{2^d(2M+1)^d}\right\rceil
\end{align}
covering $S_M\cap B_{H^\alpha}(0,1)\cap B_{L^\infty}(0,1)$ in $H^{\gamma}$. Moreover, there is an $\epsilon$-Net of size at most
\begin{align}\label{eqn:infinite_dim_eps_net}
    N_\epsilon=\left\lceil\left(\frac{4}{\epsilon}\right)^{2^d\left(2(\frac{\epsilon}{2})^{\frac{1}{\gamma}-\alpha}+1\right)^d}\right\rceil
\end{align}
covering $B_{H^\alpha}(0,1)\cap B_{L^\infty}(0,1)$ in $H^{\gamma}$.
\end{proposition}
\begin{proof} 
For a given $N_0\in\mathbb{N}$ let us denote by $\mathcal{P}_{N_0}$ the following projection
$$
\mathcal{P}_{N_0}f=\mathcal{P}_{N_0}(\sum_{k^d\in\mathbb{Z}}\hat{f}_k\e^{ikx})=\sum_{|k|\leq N_0}\hat{f}_k\e^{ikx},\text{ for any } f\in H^s(\mathbb{T}), s\geq0,
$$
where for a vector $k\in\mathbb{Z}^d$ we write \begin{align*}
|k|=\max_{l=1}^d|k_l|, \quad \langle k\rangle=\begin{cases}
    |k|,&\quad \text{if}\ k\neq0,\\
    1,&\quad\text{otherwise.}
\end{cases}    
\end{align*}
Then we have for any $u\in B_{H^\alpha}(0,1)\cap B_{L^\infty}(0,1)$:
    \begin{align*}
        \|(I-\mathcal{P}_{N_0})u\|_{H^{\gamma}}^2&=\sum_{|k|>N_0}\langle k\rangle^{2\gamma}|\hat{u}_k|^2\\
        &=\sum_{|k|>N_0}\langle k\rangle^{2\gamma}-2\alpha\langle k\rangle^{2\alpha}|\hat{u}_k|^2\leq N_0^{2(\gamma}-\alpha)\|u\|_{H^\alpha}\leq N_0^{2(\gamma}-\alpha).
    \end{align*}
    This implies that for $N_0>(\frac{\epsilon}{2})^{\frac{1}{\gamma}-\alpha}$
    \begin{align}\label{eq...firstpartcovering}
     \|(I-\mathcal{P}_{N_0})u\|_{H^{\gamma}}<\epsilon/2.
    \end{align}
    Naturally, the next step would be finding a finite $\frac{\epsilon}{2}$-net in $H^{\gamma}$ for the \textit{finite dimensional} space $\mathcal{P}_{N_0}(B_{H^\alpha}(0,1)\cap B_{L^\infty}(0,1))=S_{N_0}\cap B_{H^\alpha}(0,1)\cap B_{L^\infty}(0,1)$. Since $\|u\|_{H^{\gamma}}\leq \|u\|_{H^\alpha}$ it thus suffices to find a finite $\frac{\epsilon}{2}$-net for $\mathcal{P}_{N_0}(B_{H^{\gamma}}(0,1))$ in $H^{\gamma}$. To do so we construct an isometry that maps $\mathcal{P}_{N_0}(B_{H^{\gamma}}(0,1))$ into $\mathbb{R}^{2^d}(2N_0+1)^{d}$ and then rely on standard covering number estimates on Euclidean spaces. Our isometry of choice is the map
    \begin{align*}
        M:\mathcal{P}_{N_0}(B_{H^{\gamma}}(0,1))&\rightarrow \mathbb{C}^{(2N_0+1)^d}\\
        u&\mapsto \left(\langle k\rangle^{\gamma}\hat{u}_k\right)_{|k|\leq N_0}
    \end{align*}
    Thus it suffices to find an $\dfrac{\epsilon}{2}$-net for the ball of radius $1$ in $\mathbb{R}^{2^d}(2N_0+1)^d$.
    From the conclusion of Example 27.1 in \cite{understandML}, the smallest size of such a $\frac{\epsilon}{2}$-net, $N_{\frac{\epsilon}{2}}$, can be bounded as follows
    \begin{align}\label{eqn:finite_dim_estimate}
    N_{\frac{\epsilon}{2}}\leq \left(\frac{4}{\epsilon}\right)^{2^d(2N_0+1)^d}.
    \end{align}
    Thus, recalling \eqref{eq...firstpartcovering}, for any 
    \begin{align*}
        N_\epsilon>\left(\frac{4}{\epsilon}\right)^{2^d}(2N_0+1)^d>\left(\frac{4}{\epsilon}\right)^{2^d\left(2(\frac{\epsilon}{2})^{\frac{1}{\gamma}-\alpha}+1\right)^d}       
    \end{align*}
    we can find $u_{1},\dots, u_{N_\epsilon}\in \mathcal{P}_{N_0}(B_{H^\alpha}(0,1)\cap B_{L^\infty}(0,1))$ such that for all $u\in B_{H^\alpha}(0,1)\cap B_{L^\infty}(0,1)$ there is a $u_j\in\mathcal{P}_{N_0}(B_{H^\alpha}(0,1)\cap B_{L^\infty}(0,1))$ as the centre of those $\frac{\epsilon}{2}$-balls such that
    \begin{align*}
        \|u-u_{j}\|_{H^{\gamma}}&\leq \|\mathcal{P}_{N_0}(u-u_j)\|_{H^{\gamma}}+\|(I-\mathcal{P}_{N_0})(u-u_{j})\|_{H^{\gamma}}\\
        &\leq \|\mathcal{P}_{N_0}(u)-u_j\|_{H^{\gamma}}+\|(I-\mathcal{P}_{N_0})(u)\|_{H^{\gamma}}\\
        &< \frac{\epsilon}{2}+\frac{\epsilon}{2}\\
        &=\epsilon.
    \end{align*}
    Therefore, we can conclude that there is a $\epsilon$-net of size at most $\left\lceil\left(\frac{4}{\epsilon}\right)^{2^d\left(2(\frac{\epsilon}{2})^{\frac{1}{\gamma}-\alpha}+1\right)^d}\right\rceil$ covering $B_{H^\alpha}(0,1)\cap B_{L^\infty}(0,1)$ in $H^{\gamma}$. This concludes the proof of \eqref{eqn:infinite_dim_eps_net}. Moreover, the estimate \eqref{eqn:finite_dim_estimate} immediately implies \eqref{eqn:finite_dim_eps_net}.
\end{proof}

\section{Wall-clock estimates of individual
components in the algorithm}\label{appendix:wall-clock}
To better understand the reduction in computational time provided by the NN initialisation, we report in the following the cost of individual components of our algorithm. We highlight that the total CPU time, presented in Table \ref{table: CPUtime}, includes both the time required to compute the initial guess and the time for each iteration of Newton's method. In practice, one iteration of Newton’s method consists of three main components: the computation of the Jacobian, the calculation of the vector '\textbf{b}', and the solution of the linear system. 
Thus for the convenience of comparison, Table \ref{table:1d_wallclock} presents the time required for the offline training process,  the average time required for inferring over network, computation of Jacobian and linear solver (standard \verb|np.linalg.solve|) employed in a single step of Newton's method in 1D cases. Note that during the computation of an initial guess for the Newton iterations (online phase) training is no longer necessary and so this part does not contribute to the computational cost of the use of our algorithm in practice. Since in 1D cases, the direct assembly of the Jacobian is not particularly time-consuming, we compute the Jacobian directly in each iteration.
\begin{table}[htbp]
\tiny
\centering
\begin{tabular}{|c|c|c|c|c|c|}
\hline
Test case    & $\tau_{midpoint}$ & offline training time & inference over NN & computation of Jacobian & linear solver\\ \hline
\multirow{3}{*}{1D} & 0.5    & 2.78e+2     & 2.23e-2    & 2.15e-3   & 1.60e-2 \\ \cline{2-6} 
                    & 1     & 2.25e+2      & 1.92e-2    & 2.54e-3  &1.42e-2 \\ \cline{2-6} 
                    & 2     & 1.72e+2      & 1.61e-2    & 3.84e-3   &1.06e-2 \\ \hline
\end{tabular}
\caption{Wall-clock estimates (sec) of the training time, the average time required for inference over NN, computation of the Jacobian and linear solver in 1D cases.}
\label{table:1d_wallclock}
\end{table}

In contrast, for the 2D cases, both the computation of the Jacobian and the solution of the linear system are significantly more expensive. To improve efficiency, we avoid fully assembling the Jacobian and instead wrap the matrix-vector multiplication into a \verb|LinearOperator|, solving the linear system using the GMRES solver. This approach allows Newton's method to remain highly efficient without needing to explicitly construct the Jacobian in our implementation. Nevertheless, we report the CPU time for a single GMRES solve (used at each iteration of Newton's method) when utilizing the NN initialisation in Table \ref{table:2d_wallclock}.
\begin{table}[htbp]
\footnotesize
\centering
\begin{tabular}{|c|c|c|c|c|}
\hline
Test case           & $\tau_{midpoint}$ & offline training time & inference over NN & GMRES \\ \hline
\multirow{3}{*}{2D} & 0.25               &   4322.24             & 0.06 & 7.36   \\ \cline{2-5} 
                    & 0.5                 & 4315.04                 & 0.08 & 10.31   \\ \cline{2-5} 
                    & 1                 &      4325.51         & 0.06 &  18.11  \\ \hline
\end{tabular}
\caption{Wall-clock estimates (sec) of the training time, the average time required for inference over NN and GMRES solver in 2D cases.}
\label{table:2d_wallclock}
\end{table}

From the tables, we observe that in the 1D cases, the NN inference time is comparable to the time required for a linear solver in Newton's iteration. Therefore, if the number of iterations can be reduced by more than one, the NN initialisation leads to a notable decrease in overall computational cost. This effect is even more pronounced in the 2D cases, where the NN inference time is significantly smaller than the time required for a single GMRES solver.
 \bibliographystyle{elsarticle-num} 
 \bibliography{biblio}

\begin{thebibliography}{10}
\expandafter\ifx\csname url\endcsname\relax
  \def\url#1{\texttt{#1}}\fi
\expandafter\ifx\csname urlprefix\endcsname\relax\def\urlprefix{URL }\fi
\expandafter\ifx\csname href\endcsname\relax
  \def\href#1#2{#2} \def\path#1{#1}\fi

\bibitem{hairer1993solving}
E.~Hairer, S.~N{\o}rsett, G.~Wanner, {Solving Ordinary Differential Equations
  II: Stiff and Differential-Algebraic Problems}, Solving Ordinary Differential
  Equations II: Stiff and Differential-algebraic Problems, Springer, 1993.

\bibitem{CFL}
R.~{Courant}, K.~{Friedrichs}, H.~{Lewy}, {{\"U}ber die partiellen
  Differenzengleichungen der mathematischen Physik}, Mathematische Annalen 100
  (1928) 32--74.

\bibitem{quasinewton}
C.~G. Broyden, A class of methods for solving nonlinear simultaneous equations,
  Mathematics of Computation 19 (1965) 577--593.

\bibitem{quasireview}
E.~S. Mehiddin Al-Baali, F.~Maggioni, Broyden's quasi-newton methods for a
  nonlinear system of equations and unconstrained optimization: a review and
  open problems, Optimization Methods and Software 29~(5) (2014) 937--954.
\newblock \href {https://doi.org/10.1080/10556788.2013.856909}
  {\path{doi:10.1080/10556788.2013.856909}}.

\bibitem{inexact}
R.~S. Dembo, S.~C. Eisenstat, T.~Steihaug, Inexact newton methods, SIAM Journal
  on Numerical Analysis 19~(2) (1982) 400--408.
\newblock \href {https://doi.org/10.1137/0719025} {\path{doi:10.1137/0719025}}.

\bibitem{LIU2024112548}
L.~Liu, W.~Gao, H.~Yu, D.~E. Keyes, Overlapping multiplicative schwarz
  preconditioning for linear and nonlinear systems, Journal of Computational
  Physics 496 (2024) 112548.
\newblock \href {https://doi.org/https://doi.org/10.1016/j.jcp.2023.112548}
  {\path{doi:https://doi.org/10.1016/j.jcp.2023.112548}}.

\bibitem{chenchen}
T.~Chen, H.~Chen, Universal approximation to nonlinear operators by neural
  networks with arbitrary activation functions and its application to dynamical
  systems, IEEE Transactions on Neural Networks 6~(4) (1995) 911--917.
\newblock \href {https://doi.org/10.1109/72.392253}
  {\path{doi:10.1109/72.392253}}.

\bibitem{Lulu}
L.~Lu, P.~Jin, G.~Pang, Z.~Zhang, G.~E. Karniadakis, Learning nonlinear
  operators via {DeepONet} based on the universal approximation theorem of
  operators, Nature Machine Intelligence 3~(3) (2021) 218--229.

\bibitem{graphkernelnet}
A.~Anandkumar, K.~Azizzadenesheli, K.~Bhattacharya, N.~Kovachki, Z.~Li, B.~Liu,
  A.~Stuart, Neural operator: Graph kernel network for partial differential
  equations, in: ICLR 2020 Workshop on Integration of Deep Neural Models and
  Differential Equations, 2019.

\bibitem{li2021fourier}
Z.~Li, N.~B. Kovachki, K.~Azizzadenesheli, B.~liu, K.~Bhattacharya, A.~Stuart,
  A.~Anandkumar, Fourier neural operator for parametric partial differential
  equations, in: International Conference on Learning Representations, 2021.

\bibitem{Solver-in-the-Loop}
K.~Um, R.~Brand, Y.~Fei, P.~Holl, N.~Thuerey, Solver-in-the-loop: Learning from
  differentiable physics to interact with iterative pde-solvers, Advances in
  Neural Information Processing Systems (2020).

\bibitem{HINTS}
E.~Zhang, A.~Kahana, A.~Kopaničáková, E.~Turkel, R.~Ranade, J.~Pathak, G.~E.
  Karniadakis, \href{https://arxiv.org/abs/2208.13273}{Blending neural
  operators and relaxation methods in pde numerical solvers} (2024).
\newblock \href {http://arxiv.org/abs/2208.13273} {\path{arXiv:2208.13273}}.
\newline\urlprefix\url{https://arxiv.org/abs/2208.13273}

\bibitem{qianxiao1}
S.~Arisaka, Q.~Li, Principled acceleration of iterative numerical methods using
  machine learning, in: Proceedings of the 40th International Conference on
  Machine Learning, ICML, 2023.

\bibitem{qianxiao2}
S.~Arisaka, Q.~Li, Accelerating legacy numerical solvers by non-intrusive
  gradient-based meta-solving, in: Proceedings of the 40th International
  Conference on Machine Learning, ICML, 2024.

\bibitem{IntDeep}
J.~Huang, H.~Wang, H.~Yang, Int-deep: A deep learning initialized iterative
  method for nonlinear problems, Journal of Computational Physics 419 (2020)
  109675.
\newblock \href {https://doi.org/https://doi.org/10.1016/j.jcp.2020.109675}
  {\path{doi:https://doi.org/10.1016/j.jcp.2020.109675}}.

\bibitem{doi:10.1137/22M1507942}
L.~Luo, X.-C. Cai, \(\text{PIN}^{\mathcal l}\) : Preconditioned inexact newton
  with learning capability for nonlinear system of equations, SIAM Journal on
  Scientific Computing 45~(2) (2023) A849--A871.
\newblock \href {https://doi.org/10.1137/22M1507942}
  {\path{doi:10.1137/22M1507942}}.

\bibitem{co2}
A.~Lechevallier, S.~Desroziers, T.~Faney, E.~Flauraud, F.~Nataf, {Hybrid Newton
  method for the acceleration of well event handling in the simulation of CO2
  storage using supervised learning}, working paper or preprint (Apr. 2023).

\bibitem{aghili2024accelerating}
J.~Aghili, E.~Franck, R.~Hild, V.~Michel-Dansac, V.~Vigon, Accelerating the
  convergence of newton's method for nonlinear elliptic pdes using fourier
  neural operators (2024).
\newblock \href {http://arxiv.org/abs/2403.03021} {\path{arXiv:2403.03021}}.

\bibitem{tierra_review_2015}
G.~Tierra, F.~Guill{\'e}n-Gonz{\'a}lez, {Numerical Methods for Solving the
  Cahn--Hilliard Equation and Its Applicability to Related Energy-Based
  Models}, Archives of Computational Methods in Engineering 22~(2) (2015)
  269--289.

\bibitem{timestepadapative}
F.~Guillén-González, G.~Tierra, {Second order schemes and time-step
  adaptivity for Allen–Cahn and Cahn–Hilliard models}, Computers \&
  Mathematics with Applications 68~(8) (2014) 821--846.
\newblock \href {https://doi.org/https://doi.org/10.1016/j.camwa.2014.07.014}
  {\path{doi:https://doi.org/10.1016/j.camwa.2014.07.014}}.

\bibitem{gomez2011}
H.~Gomez, T.~J. Hughes, {Provably unconditionally stable, second-order
  time-accurate, mixed variational methods for phase-field models}, Journal of
  Computational Physics 230~(13) (2011) 5310--5327.

\bibitem{hou2023energy}
D.~Hou, L.~Ju, Z.~Qiao, Energy-dissipative spectral renormalization exponential
  integrator method for gradient flow problems, arXiv preprint arXiv:2310.00824
  (2023).

\bibitem{AC_im_midpoint}
A.~Rodgers, D.~Venturi, Implicit integration of nonlinear evolution equations
  on tensor manifolds, Journal of Scientific Computing 97~(2) (2023) 33.

\bibitem{mishra22}
S.~Mishra, R.~Molinaro, {Estimates on the generalization error of
  physics-informed neural networks for approximating PDEs}, IMA Journal of
  Numerical Analysis 43~(1) (2022) 1--43.

\bibitem{BoundaryConditionsCNN}
A.~Alguacil, W.~G. Pinto, M.~Bauerheim, M.~C. Jacob, S.~Moreau, Effects of
  boundary conditions in fully convolutional networks for learning
  spatio-temporal dynamics, in: Machine Learning and Knowledge Discovery in
  Databases. Applied Data Science Track, Springer International Publishing,
  Cham, 2021, pp. 102--117.

\bibitem{CNNreview}
X.~Zhao, L.~Wang, Y.~Zhang, X.~Han, M.~Deveci, M.~Parmar, A review of
  convolutional neural networks in computer vision, Artificial Intelligence
  Review 57~(4) (2024) 99.

\bibitem{pytorch}
A.~Paszke, S.~Gross, F.~Massa, A.~Lerer, J.~Bradbury, G.~Chanan, T.~Killeen,
  Z.~Lin, N.~Gimelshein, L.~Antiga, A.~Desmaison, A.~K\"{o}pf, E.~Yang,
  Z.~DeVito, M.~Raison, A.~Tejani, S.~Chilamkurthy, B.~Steiner, L.~Fang,
  J.~Bai, S.~Chintala, PyTorch: an imperative style, high-performance deep
  learning library, Curran Associates Inc., 2019.

\bibitem{adam}
D.~P. Kingma, J.~Ba, Adam: {A} method for stochastic optimization, in: 3rd
  International Conference on Learning Representations, 2015.

\bibitem{SuliMayers2003}
E.~Süli, D.~F. Mayers, An Introduction to Numerical Analysis, Cambridge
  University Press, 2003.

\bibitem{1975}
R.~A. Adams, J.~J. Fournier, Sobolev Spaces, Academic Press, 1975.

\bibitem{finite_dim_appro_error}
K.~Bhattacharya, B.~Hosseini, N.~B. Kovachki, A.~M. Stuart, Model {Reduction}
  {And} {Neural} {Networks} {For} {Parametric} {PDEs}, The SMAI Journal of
  computational mathematics 7 (2021) 121--157.
\newblock \href {https://doi.org/10.5802/smai-jcm.74}
  {\path{doi:10.5802/smai-jcm.74}}.

\bibitem{DeepOnet_error}
S.~Lanthaler, S.~Mishra, G.~E. Karniadakis, {Error estimates for DeepONets: a
  deep learning framework in infinite dimensions}, Transactions of Mathematics
  and Its Applications 6~(1) (2022) tnac001.
\newblock \href {https://doi.org/10.1093/imatrm/tnac001}
  {\path{doi:10.1093/imatrm/tnac001}}.

\bibitem{FNO_error}
T.~Kim, M.~Kang, Bounding the rademacher complexity of fourier neural
  operators, Machine Learning 113~(5) (2024) 2467--2498.

\bibitem{CMS_weinan}
W.~E, C.~Ma, Q.~Wang, Rademacher complexity and the generalization error of
  residual networks, Communications in Mathematical Sciences 18 (2020)
  1755--1774.
\newblock \href {https://doi.org/10.4310/CMS.2020.v18.n6.a10}
  {\path{doi:10.4310/CMS.2020.v18.n6.a10}}.

\bibitem{DBLP:journals/corr/SzegedyZSBEGF13}
C.~Szegedy, W.~Zaremba, I.~Sutskever, J.~Bruna, D.~Erhan, I.~J. Goodfellow,
  R.~Fergus, Intriguing properties of neural networks, in: International
  Conference on Learning Representations, 2014.

\bibitem{ALLEN19791085}
S.~M. Allen, J.~W. Cahn, A microscopic theory for antiphase boundary motion and
  its application to antiphase domain coarsening, Acta Metallurgica 27~(6)
  (1979) 1085--1095.
\newblock \href {https://doi.org/https://doi.org/10.1016/0001-6160(79)90196-2}
  {\path{doi:https://doi.org/10.1016/0001-6160(79)90196-2}}.

\bibitem{im-exp}
T.~Tang, J.~Yang, Implicit-explicit scheme for the allen-cahn equation
  preserves the maximum principle, Journal of Computational Mathematics 34~(5)
  (2016) 451--461.

\bibitem{doi:10.1137/19M1243750}
Q.~Du, L.~Ju, X.~Li, Z.~Qiao, Maximum bound principles for a class of
  semilinear parabolic equations and exponential time-differencing schemes,
  SIAM Review 63~(2) (2021) 317--359.
\newblock \href {https://doi.org/10.1137/19M1243750}
  {\path{doi:10.1137/19M1243750}}.

\bibitem{Goeckler2014}
T.~G{\"{o}}ckler, {Rational Krylov subspace methods for phi-functions in
  exponential integrators}, Ph.D. thesis, {Karlsruher Institut f{\"{u}}r
  Technologie (KIT)} (2014).
\newblock \href {https://doi.org/10.5445/IR/1000043647}
  {\path{doi:10.5445/IR/1000043647}}.

\bibitem{understandML}
S.~Shalev-Shwartz, S.~Ben-David, Understanding Machine Learning - From Theory
  to Algorithms., Cambridge University Press, 2014.

\end{thebibliography}





\end{document}